\documentclass[dvips,twoside]{article}
\usepackage[latin1]{inputenc}
\usepackage[T1]{fontenc}
\usepackage{parskip}
\usepackage{amsmath,amsfonts,amssymb,amsxtra}
\usepackage{latexsym}
\usepackage{theorem}
\usepackage{graphicx,psfrag}
\usepackage{verbatim}
\usepackage{hyperref}

{\theoremstyle{plain}                       
\theorembodyfont{\itshape}
\newtheorem{Theorem}{Theorem}}

{\theoremstyle{plain}                       
\theorembodyfont{\itshape}
\newtheorem{Corollary}{Corollary}}

{\theoremstyle{plain}                       
\theorembodyfont{\itshape}
\newtheorem{Proposition}{Proposition}}

{\theoremstyle{plain}                       
\theorembodyfont{\rmfamily}
\newtheorem{Definition}{Definition}}

{\theoremstyle{plain}                       
\theorembodyfont{\rmfamily}
}

{\theoremstyle{plain}                       
\theorembodyfont{\rmfamily}
\newtheorem{Example}{Example}}

{\theoremstyle{plain}                       
\theorembodyfont{\rmfamily}
\newtheorem{Remark}{Remark}}

{\theoremstyle{plain}                       
\theorembodyfont{\itshape}
\newtheorem{Lemma}{Lemma}}

{\theoremstyle{plain}                       
\theorembodyfont{\itshape}
\newtheorem{Condition}{Condition}}

\def\M{\mathcal{M}}
\def\N{\mathbb{N}}
\def\Z{\mathbb{Z}}

\begin{document}

\title{Equilibrium states and invariant measures  for  random dynamical systems}
\author{Ivan Werner\\
   {\small Email: ivan\_werner@mail.ru}}
\date{November 9, 2014}
\maketitle

\begin{abstract}\noindent
Random dynamical systems with countably many maps which admit countable Markov partitions  on complete metric spaces such that the resulting Markov systems are uniformly continuous and contractive are considered. A non-degeneracy and a consistency conditions for such systems, which admit some proper Markov partitions of connected spaces, are introduced, and further sufficient conditions for them are provided.  It is shown that every uniformly continuous Markov system associated with a continuous random dynamical system  is consistent if it has a dominating Markov chain. A  necessary and sufficient condition for the existence of an invariant Borel probability measure for such a non-degenerate system with a  dominating Markov chain and a finite \eqref{cgc} is given.  The condition is also sufficient if the non-degeneracy is weakened with the consistency condition. A further sufficient condition for the existence of an invariant measure for such a consistent system which involves only the properties of the  dominating Markov chain  is provided. In particular, it implies that every such  a consistent system with a finite Markov partition and a finite \eqref{cgc} has an invariant Borel probability measure. A bijective map between these measures and equilibrium states associated with such a  system is established in the non-degenerate case. Some properties of the map and the measures are given.


\end{abstract}

\tableofcontents

\section{Introduction}
The purpose of this note is to show the existence of invariant measures for some random dynamical systems introduced in \cite{Wer1} as {\it contractive Markov systems}. To avoid some confusion, let us stress that the word 'Markov' in the name was used to indicate a Markovian topological structure of the random dynamical system, which naturally generalizes a weighted directed graph, but the dependence structure of random processes which can be generated by it on a code space, in contrast to directed graphs,  can be far beyond Markovian. They certainly can generate any stationary  process with values in a discrete state space. At the same time, the algorithm for generating a process by such a system is not much different to that for generating a process by a weighted directed graph. Such random dynamical systems find more and more applications in  modern sciences, e.g. \cite{Sl}, \cite{BLLC}, \cite{FHK}, and provide new challenging and illuminating examples for mathematical theories, e.g. \cite{Wer6}, \cite{Wer10}, \cite{Wer12}. However, in contrast to  weighted  directed graphs, the behaviour of contractive Markov systems is still not fully understood. 

The existence of stationary states for such systems was shown  on some locally compact spaces in  \cite{Wer1}. This  was proved  under the condition that the partition of the Markov system consists of open sets. Though, this was sufficient to cover finite {\it Markov chains} and $g$-{\it measures} \cite{Ke} with the theory, it clearly poses a severe restriction on  the applicability of it. In particular, the removal of the condition admits the usage of {\it Markov partitions} for  random dynamical systems, which reduces the latter to Markov systems, the behaviour of which is more transparent \cite{Wer11}. The proof which was given in \cite{Wer1} went along the lines of that which had been given by M. Barnsley et al. for {\it iterated function systems with place-dependent probabilities}  \cite{BDEG} and \cite{BDEGE}. 
The result then was extended  by K. Horbacz and T. Szarek  \cite{HS}  on Polish spaces, through  application of some results  which had been obtained by the second author  for Markov operators satisfying some {\it non-expansiveness
} and {\it concentration} conditions on Polish spaces, using a {\it lower bound technique}  \cite{Sz}. Unfortunately, the condition of the openness of the partition has been left in place. 

In this article, we close the gap by introducing a {\it non-degeneracy}  and a {\it consistency} conditions and providing further sufficient conditions for  them. These conditions admit some proper Markov partitions of connected spaces and allow us to prove the existence of invariant measures for the random dynamical systems which exhibit the continuity and the contraction on average properties only on the atoms of their Markov partitions. In particular, the consistency condition includes all uniformly continuous Markov systems which  are associated with  random dynamical systems with continuous maps and probability functions and have a {\it dominating Markov chain} (Theorem \ref{ucct}).  Furthermore, every countable refinement of a uniformly continuous, positive and consistent Markov system with contractive maps is again consistent if it has a dominating Markov chain (Proposition \ref{crp}).
 
The presented proof is self-contained, does not require any special knowledge and works for countable Markov systems. Moreover, it is shown that the separability of the space is not needed in this case. The existence of the invariant measures is deduced from the existence of {\it equilibrium states} or, in general,  {\it asymptotic states} on the code space associated with such a system, via a {\it coding map}. This method is easier because the code space is either a compact metrizable space, as in  the case with finitely many maps, or can be easily extended to such a space in the case of countably many maps. In particular, one can take advantage of the weak-star compactness of the set of all Borel probability measures on it. 

The  existence of the equilibrium states for  {\it energy functions} associated with such systems  has been  already shown in \cite{Wer6}, but it has been deduced from the existence of the invariant measures for such systems with the open partition on locally compact spaces. The main message  of \cite{Wer6} was  that the current thermodynamic formalism is not applicable to such systems because it fails even to predict  the existence of  equilibrium states for such energy functions, not to mention the construction of them. 

Recently, a construction of such equilibrium states has been  proposed in \cite{Wer10} and \cite{Wer12}. It seems to require the existence of an equilibrium state for such a system for the proof that the constructed measure is not zero. 

An other related existence result is a recent proof in the particular case of $g$-measures by A. Johansson et al. \cite{JOP}. However, it does not intersect much with the present result as the $g$-functions associated with  our systems are not continuous, even in the case of the openness of the Markov partition (see \cite{Wer6}) or the case of contractive maps on a compact metric space, but  without the openness condition on the Markov partition (e.g. see Example \ref{mcae} below).

The present result also establishes a bijection between the equilibrium states and the invariant Borel probability measures of such systems in the non-degenerate case. In particular, this generalizes a theorem by F. Ledrappier \cite{Le},  Theorem 2.1 in \cite{Wal}.

The article is organized as follows. Section 1 collects all the necessary definitions and notations. Section 2 presents the main results. Finally, Section 3 provides, in particular, some simple examples  to which the technique of  Markov partitions can be applied. As far as the author is aware, some of the examples have not been accessible by the theory before. 

It was pointed out by an anonymous reviewer that it might be appropriate to cite the works \cite{MU},\cite{SU},\cite{MSU},\cite{RU}, despite the fact that the authors of the works seem to indicate that their works  are not related to contractive Markov systems, by not citing any previous works on contractive Markov systems. It is a remarkable coincidence that the 'graph directed constructions' \cite{MW} evolved into the 'graph directed Markov systems' \cite{MU} in literature, even on the costs of making the new name somewhat tautological, shortly after the author had introduced the 'contractive Markov systems' with contractive maps on compact metric spaces in his diploma thesis \cite{Wer0} (certainly, aware of \cite{MW}, but mostly influenced  by the work of J. Elton \cite{Elton}, see \cite{Wer11}). (As a matter of fact, the diploma thesis was then developed further to his Ph.D. thesis at the University of St Andrews in 2003, but for some reason the university allowed to defend the thesis only after 1 year of waiting, in November 2004.) Notable also is that the unusual direction of the arrows of the directed graphs in \cite{MW} evolved also into that of Markov systems. 
(Curiously, the anonymous reviewer pointed out that a term 'conformal graph directed Markov systems' had already  appeared in \cite{MU0}, but apparently only in the first two sentences as a name for a future theory towards which the authors direct their efforts without giving a definition for such objects yet. In fact,  the rest of the article seems to be completely detached from the information contained in these two sentences. However, the date of that publication and the date when the author was allowed to officially submit the diploma thesis make the coincidence even more remarkable.)
The author leaves it to the reader to judge on how the structures studied in the cited works relate to contractive Markov systems  and on their scientific motivation. 

For some criteria for the uniqueness of the invariant probability measures, the existence of which is proved in this article, the reader is referred to \cite{Wer13}. 

\section{Definitions and notation}

 Let $\mathcal{B}(X)$ denote the Borel $\sigma$-algebra on a topological space $X$  and $P(X)$ denote the set of all Borel probability measures on it. Let $\mathcal{L}^B(X)$ denote the set of all real-valued, non-negative, Borel measurable functions on $X$. For $B\in\mathcal{B}(X)$, let $P(B)$ denote the set of all $\nu\in P(X)$ such that $\nu(B)=1$.

Let $(K,d)$ be a complete metric space. A family $D_R:=(K,w_e,p_e)_{e\in E}$ is called a  {\it random dynamical system} on $K$ iff $E$ is at most countable,  $w_e:K\longrightarrow K$ and $p_e:K\longrightarrow [0,1]$ are Borel-measurable for all $e\in E$ such that $\sum_{e\in E}p_e(x) = 1$ for all $x\in K$. $w_e$'s are called {\it maps} and $p_e$'s are called {\it probability functions}. 
 
With $D_R$ is associated a {\it Markov operator} $U$ defined on $\mathcal{L}^B(K)$ by 
\[Uf:=\sum\limits_{e\in E}p_ef\circ w_e\]
for all $f\in\mathcal{L}^B(K)$. Let $U^*$ denote its {\it adjoint operator} acting on $\nu\in P(K)$ by $U^*\nu(f):=\int Ufd\nu$ for all  bounded $f\in \mathcal{L}^B(K)$. $\mu\in P(K)$ is called an {\it invariant measure} for the random dynamical system iff $U^*\mu=\mu$. Observe that each $w_e$ needs to be defined only on the set $\{x\in K|\ p_e(x)> 0\}$ for the definitions of $U$ and $U^*$, it can be then extended on the whole space arbitrarily.

A random dynamical system is called a  {\it  Markov system}  if and only if it has the form $(K_{i(e)},w_e, p_e)_{e\in E'}$  where $E'$ is a  set such that there exists a partition of $K$ into non-empty Borel subsets  $(K_j)_{j\in N}$, $N\subset\mathbb{N}$ with $1\in N$ (case where the size of $ N$ is $1$ is not excluded),  and a surjective map $i:E'\longrightarrow N$  and $t:E'\longrightarrow N$ such that for every  $e\in E'$  there exist Borel measurable $w_e:K_{i(e)}\longrightarrow K_{t(e)}$ and $p_e:K_{i(e)}\longrightarrow [0, 1]$  such that there exists $x_e\in K_{i(e)}$ with $p_e(x_e) >0$, and $\sum_{e\in E', i(e)=j}p_e(y) = 1$ for all $y\in K_j$ and $j\in N$. $K_j$'s are called the {\it vertex sets} of the Markov system. The Markov system is called {\it countable} iff  $N$ and $E$ are at most countable. Clearly, a countable Markov system defines a random dynamical system on $K$ by extending $p_e$'s on $K$ by zero and $w_e$'s arbitrarily. Such extensions define the actions of the Markov system on functions and measures through operators $U$ and $U^*$ and will be always  assumed.

We say that a random dynamical system has a {\it Markov partition} iff there exists  a partition of $K$ into non-empty Borel subsets  such that the restrictions of its maps and probability functions on the atoms of the partition (after a possible re-indexation) form a Markov system.

A Markov system  $(K_{i(e)},w_e, p_e)_{e\in E}$ is called {\it contractive} with a contraction rate $0<a<1$ iff  
 \begin{equation}\label{cac}
 \sum\limits_{e\in E, i(e)=j}p_e(x)d(w_ex,w_ey)<ad(x,y)\\\\\mbox{ for all }x,y\in K_j\mbox{ and } j\in N.
 \end{equation}
The condition was introduced by R. Isaac in \cite{Is} for the case of $N=\{1\}$.

We say that a Markov system  $(K_{i(e)},w_e, p_e)_{e\in E}$ is {\it (uniformly) continuous} iff  $w_e|_{K_{i(e)}}$ and  $p_e|_{K_{i(e)}}$ are (uniformly) continuous for each $e\in E$, where the notation $f|_A$ means the restriction of a function $f$ on a set $A$. We call the Markov system {\it positive} iff $p_e|_{K_{i(e)}}>0$ for all $e\in E$.

A sequence $(e_1,...,e_n)$ of $e_i\in E$ for all $1\leq i\leq n$ is called a {\it path} of the Markov system iff $i(e_{i+1})=t(e_i)$ for all $i$.
We will denote by $\delta_x\in P(K)$  the Dirac probability measure concentrated at $x\in K$, by $B_\alpha(x)$ the closed ball of  radius $\alpha$ and centre $x$, by  $1_A$  the indicator function of a set $A$, by $\bar A$ the topological closure of a set $A$ and by $\bar f$ the continuous extension of a uniformly continuous function $f$ on the closure of the domain of its definition. For a measurable map between measure spaces
$f:(X,\mathcal{A}, \mu)\longrightarrow (Y,\mathcal{B})$, $f(\mu)$ will denote the measure on $ (Y,\mathcal{B})$ given by $f(\mu)(B):=\mu(f^{-1}(B))$ for all $B\in\mathcal{B}$, and $f^{-1}(\mathcal{B})$ will denote the $\sigma$-algebra   $\{f^{-1}(B)|\ B\in\mathcal{B}\}$. As usual, $\ll$ will denote the absolute continuity relation for measures.

Let $(X, \mathcal{B}, \Lambda)$ be a probability space and $I$ be an at most countable set.  A family $(A_i)_{i\in I}\subset\mathcal{B}$ is called a {\it partition} of $(X, \mathcal{B}, \Lambda)$ iff its members are pairwise disjoint and $\Lambda(X\setminus\bigcup_{i\in I}A_i)=0$.  For a partition $\alpha$ of  $(X, \mathcal{B}, \Lambda)$ and a sub-$\sigma$-algebra $\mathcal{C}\subset\mathcal{B}$, $H_\Lambda(\alpha|\mathcal{C})$ will denote the {\it conditional entropy} of $\alpha$ conditioned on $\mathcal{C}$ with respect to $\Lambda$, which is given by
\[H_\Lambda(\alpha|\mathcal{C}):=-\sum\limits_{A\in \alpha}\int E_\Lambda(1_{A}|\mathcal{C})\log E_\Lambda(1_{A}|\mathcal{C})d\Lambda,\]
with the usual definition $0\log 0:=0$, where $E_\Lambda(1_{A}|\mathcal{C})$ denotes the conditional expectation of the indicator function $1_A$  conditioned on $\mathcal{C}$ with respect to  $\Lambda$.

We will use the usual notion of the tightness. A set of Borel measures $\{\Lambda_i|\ i\in I\}$ on a topological space $X$ is called {\it(uniformly) tight } iff for every $\epsilon>0$ there exists a compact $C\subset X$ such that $\Lambda_i(X\setminus C)<\epsilon$ for all $i\in I$.

\section{Results}

Let $D_R$ be a random dynamical system on a complete metric space $(K, d)$ which has a Markov partition $(K_j)_{j\in N}$ such that the resulting Markov system $\mathcal{M}:=(K_{i(e)},w_e, p_e)_{e\in E}$ is countable. Set \[P(\mathcal{M}):=\left\{\mu\in P(K)|\ U^*\mu=\mu \right\}.\]

Let $E$ and $N$ be provided with the discrete topologies. Set $\bar E:=E\cup\{\infty\}$ endowed with Alexandrov's one-point compactification topology, i.e.  the topology consists of all subsets of $E$ and sets of the form $\bar E\setminus C$ where  $C\subset E$ is finite. Note that the topology has a countable base (the axiom of choice is assumed in this paper). Let $\bar E$ be  equipped with the Borel $\sigma$-algebra. Note that the Borel $\sigma$-algebra still consist of all subsets. We can write $D_R=(K, w_e, p_e)_{e\in\bar E}$ where $w_\infty := id$ and $p_\infty := 0$, as such an extension does not change  the action of $D_R$ on functions and measures by its operators.  Set $i(\infty) := 1$ and $t(\infty):= 1$. Then we also can write $\mathcal{M}=(K_{i(e)},w_e, p_e)_{e\in\bar E}$ in the above sense. Now, set $\bar\Sigma:=\{\sigma:=(...,\sigma_{-1},\sigma_{0},\sigma_{1},...)|\ \sigma_i\in \bar E\mbox{ for all }i\in\mathbb{Z}\}$ provided with the product topology (a similar compactification has been used in \cite{JOP}). Hence $\bar\Sigma$ is Hausdorff and, by the Tikhonov Theorem, compact. Moreover, the topology of $\bar\Sigma$ has a countable base, since $\bar E$ does, and it is regular, since $\bar E$ is, and therefore, it is metrizable, by the Urysohn Metrization Theorem.  $\bar\Sigma$ is called the {\it code space} of the Markov system. Note that, since the topology of $\bar E$ has a countable base, the Borel $\sigma$-algebra on $\bar\Sigma$ coincides with the product $\sigma$-algebra.  Let $S:\bar\Sigma\longrightarrow\bar\Sigma$ be the {\it left shift} map given  by $(S\sigma)_{i-1} = \sigma_i$ for all $i\in\mathbb{Z}$ and $\sigma\in\bar\Sigma$. Let $P_S(\bar\Sigma)$ denote the space of all shift invariant Borel probability measures on $\bar\Sigma$ equipped with the weak-star topology. Recall that, since  the topology of $\bar\Sigma$ has a countable base, the Banach space of all continuous functions on it is separable, and therefore, the weak-star topology on the unit ball of the dual space is  metrizable. Furthermore, by the Riesz Representation Theorem and the Alaoglu Theorem,  $P_S(\bar\Sigma)$ is compact and metrizable in the weak-star topology, as a closed subset of the unit ball.

Let $m\leq n\in\mathbb{Z}$ and $e_m,...,e_n\in\bar E$. Set $_m[e_m,...,e_n]:=\{\sigma\in\bar\Sigma|\ \sigma_i=e_i\mbox{ for all }m\leq i\leq n\}$, it is called a {\it cylinder set}. Let $\mathcal{A}_m$ denote the $\sigma$-algebra generated by the cylinder sets of the form $_m[e_m,...,e_n]$, $n\geq m$, and $\mathcal{F}_m\subset\mathcal{A}_m$, $m\leq 0$, denote the  $\sigma$-algebra generated by the cylinder sets of the form $_m[e_m,...,e_0]$. Let $\mathcal{F}$ denote the $\sigma$-algebra generated by $\bigcup_{m\leq 0}\mathcal{F}_m$.

For  $x\in K$, let $P^m_x$ denote the probability measure on $\mathcal{A}_m$ given by
\[P^m_x( _m[e_m,...,e_n]):=p_{e_m}(x)p_{e_{m+1}}(w_{e_m}x)...p_{e_n}(w_{e_{n-1}}\circ ...\circ w_{e_{m}}x)\]
for all $_m[e_m,...,e_n]\in\mathcal{A}_m$, $n\geq m$, e.g. by the Kolmogorov Consistency Theorem. Observe that $P^m_x=P^0_x\circ S^{m}$ for all $m\leq 0$, $x\in K$ (here, $S^m$ denotes  the naturally induced set map).  

Let $\bar\Sigma^+:=\{(\sigma_1, \sigma_2,...)|\ \sigma_i\in \bar E\mbox{ for all }i\in\mathbb{N}\}$  provided with the product topology and the product $\sigma$-algebra, and let $\mathcal{B}(K)\otimes\mathcal{B}(\bar\Sigma^+)$ denote the product $\sigma$-algebra of the Borel $\sigma$-algebra on $K$ and that on $\bar\Sigma^+$. Let $_m[e_m,...,e_n]^+\subset\bar\Sigma^+$, $m>0$, denote a cylinder set. For $x\in K$, let $P_x$ denote the Borel probability measure on $\bar\Sigma^+$ given by 
\[P_x( _1[e_1,...,e_n]^+):=p_{e_1}(x)p_{e_{2}}(w_{e_1}x)...p_{e_n}(w_{e_{n-1}}\circ ...\circ w_{e_{1}}x)\]
for all $_1[e_1,...,e_n]^+\subset\bar\Sigma^+$.

Set
\[\Sigma_G :=\left\{\sigma\in\bar\Sigma |\  i(\sigma_{n+1}) = t(\sigma_n),\ \sigma_n\in E\mbox{ for all }n\in\mathbb{Z}\right\}\]
provided with the metric $d'(\sigma, \sigma'):=2^{-k}$ where $k\in\mathbb{N}\cup\{0\}$ is the largest with $\sigma_i=\sigma'_i$ for all $|i|<k$ for all $\sigma,\sigma'\in\Sigma_G$. We call $\Sigma_G$ the {\it path space} associate with $\M$.
One easily checks that the  topology on $\Sigma_G$ which is induced from $\bar\Sigma$ coincides with that given by  $d'$.

\begin{Remark}\label{ppr}
In the following, often implicitly, the following fact will be used, which might be useful to observe before. For every $_m[e_m,...,e_n]\in\mathcal{A}_m$, $x\in K$ and $m\in\mathbb{Z}$, $P^m_x(  _m[e_m,...,e_n])>0$ implies that $(e_m,...,e_n)$ is a {\it path} of the Markov system and $x\in K_{i(e_m)}$. This follows from the definition that the probability functions are zero outside their vertex sets.  Note that, in this paper, they are allowed to take the value zero also on their vertex sets, whereas in \cite{Wer1} it was required that $p_{e}|_{K_{i(e)}}>0$ for all $e\in E$. The latter is necessary if one wants to prove that the process started at any $x\in  K_i$, for a fixed $i\in N$, converges to the same stationary state \cite{Wer11}. However, in this article, we are concerned only with the question on the existence of the stationary states.  

Moreover, observe that, since $i$ is surjective,
   \begin{equation}\label{isur}
     \sum\limits_{(e_m,...,e_n)\mbox{ is a path}}P^{m-k}_x(  _m[e_m,...,e_n]) = 1
   \end{equation}
for all $x\in K$, $m<n$ and $k\geq 0$.
\end{Remark}

\subsection{Refinement of a Markov system}

The study of a random dynamical system via an associated Markov system has some flexibility. Often one can choose several Markov systems associated with a given random dynamical system, e.g. see Examples \ref{prdm} and \ref{dMse}. In such a case, choosing a finer Markov partition can help to obtain a Markov system with desired properties. Also, sometimes one can obtain some proprieties of a Markov system from those which refine it, e.g. see Example 1 in \cite{Wer13}, and vice versa, e.g. see Proposition \ref{crp}. In this subsection, we provide some tools which enable one to exploit this flexibility.

\begin{Definition}
    We call a Markov system  $\M^r:=(K^r_{i(e)},w^r_e,p^r_e)_{e\in E^r}$ a {\it refinement} of $\M$ if and only if partition 
$\{K^r_{i(e)}\}_{e\in E^r}$ refines  partition $\{K_{i(e)}\}_{e\in E}$ (i.e. each $K_i$ is a union of some $K^r_j$'s) and there is a surjective map $r:E^r\longrightarrow E$ such that
$w_{r(e)}|_{K^r_{i(e)}} = w^r_{e}|_{K^r_{i(e)}}$ and $p_{r(e)}|_{K^r_{i(e)}} = p^r_{e}|_{K^r_{i(e)}}$ for all $e\in E^r$ (we use the same letters for maps $i,t:E^r\longrightarrow N^r$). We call $r$  the {\it refinement map}. Let $r$ be extended on the one-point compactification by $r(\infty):=\infty$. Then $r$ defines a Borel-Borel-measurable surjective map
 \begin{eqnarray*}
       \Psi_r:&{\bar\Sigma}^r&\longrightarrow\bar\Sigma\\
                 &(\sigma_k)_{k\in\Z}&\longmapsto (r(\sigma_k))_{k\in\Z}
 \end{eqnarray*}
where ${\bar\Sigma}^r$ denotes the compact code space associated with $\M^r$, and in the same way $\psi_r:{\bar\Sigma^{r+}}\longrightarrow \bar\Sigma^+$. We will denote most objects associated with $\M^r$ with the same letters as for those associated with $\M$ and a superscript $^r$ or a subscript $_r$ (e.g. $\Sigma^r_G$ denotes the path space of  $\M^r$). 
\end{Definition}
Note that a measure is invariant for $\M$ if and only if it is invariant for $\M^r$, as with both Markov systems is associated the same Markov operator $U$.

\begin{Lemma}\label{prm}
   Suppose $\M^r$ is countable. \\
  (i) $\Psi_r$ is continuous if and only if $r^{-1}\{e\}$ is finite for all $e\in E$.\\
  (ii) $\Psi_r\left(\Sigma^r_G\right)\subset\Sigma_G$.\\
  (iii) $\Psi_r\left(\Sigma^r_G\right)=\Sigma_G$ if $\M$ is positive and $\Psi_r$ is continuous.\\
  (iv) $P_x=\psi_r(P^r_x)$ for all $x\in K$.\\
   (v)  $\Phi(\mu)=\Psi_r(\Phi_r(\mu))$ for all $\mu\in P(\M)$.
\end{Lemma}
{\it Proof.}  (i) For the 'if',  it is sufficient to show that $\Psi_r^{-1}(U)$ is open for every $U$ from the subbase of the topology on $\bar\Sigma$. 

Let $U=\L_i[e]$  for some $e\in E$ and $i\in\Z$. Then, by the definition of $\Psi_r$, $\Psi_r^{-1}(U)=\bigcup_{e'\in E^r, r(e')=e}\L_i[e']$, and therefore, it is open.

Now, let $U=\bigcup_{e\in\bar E\setminus C}\L_i[e]$ for some finite $C\subset E$ and $i\in\Z$. Then
 \begin{eqnarray*}
      \Psi_r^{-1}(U)&=&\bigcup\limits_{e\in\bar E\setminus C} \Psi_r^{-1}(\L_i[e])=\left\{\sigma'\in\bar\Sigma^r|\ r(\sigma'_i)\in\bar E\setminus C\right\}\\
&=&\left\{\sigma'\in\bar\Sigma^r|\ \sigma'_i\in\bar {E^r}\setminus r^{-1}(C)\right\}=\bigcup\limits_{e'\in\bar {E^r}\setminus r^{-1}(C)}\L_i[e'].
 \end{eqnarray*}
Hence, by the hypothesis, $ \Psi_r^{-1}(U)$ is open.

For the 'only if', suppose $\Psi_r$  is continuous. Then $\Psi_r^{-1}(U)$ is open for any of the above cases, but this is possible only if $r^{-1}(C)$ is finite (since there exists finite $C'\subset E^r$ such that $\bar {E^r}\setminus C'\subset \bar {E^r}\setminus r^{-1}(C)$, i.e. $r^{-1}(C)\subset C'$).

(ii) Let $\sigma'\in\Sigma^r_G$ and $\sigma:=\Psi_r(\sigma')$. Let $n\in\Z$ and $j:=i(\sigma_n)$. Observe that
\[w_{\sigma_n}\left(K^r_{i(\sigma'_n)}\right)=w^r_{\sigma'_n}\left(K^r_{i(\sigma'_n)}\right)\subset K^r_{t(\sigma'_n)}= K^r_{i(\sigma'_{n+1})}\subset K_{i(\sigma_{n+1})}\]
and
\[w_{\sigma_n}\left(K^r_{i(\sigma'_n)}\right)\subset w_{\sigma_n}\left(K_{i(\sigma_n)}\right)\subset K_{t(\sigma_{n})}.\]
Hence $i(\sigma_{n+1})= t(\sigma_{n})$. The assertion  follows.

(iii) Let $\sigma\in\Sigma_G$. Since $\Psi_r$ is surjective, there exists $\sigma'\in\bar\Sigma^r$ such that $\Psi_r(\sigma')=\sigma$. Suppose there exists $j\in\Z$ such that $t(\sigma'_j)\ne i(\sigma'_{j+1})$. Choose such $j$ with the smallest absolute value. In the case that there are two such $j$, do what follows first for the non-positive one and then iterate the choice of such $j$ with the smallest absolute value again. By the definition of $\M^r$,
\[w^r_{\sigma'_j}\left(K^r_{i(\sigma'_j)}\right)=w_{\sigma_j}\left(K^r_{i(\sigma'_j)}\right)\subset K_{t(\sigma_j)}=K_{i(\sigma_{j+1})}.\]
Hence
\[w^r_{\sigma'_j}\left(K^r_{i(\sigma'_j)}\right)\subset K_{i(\sigma_{j+1})}\setminus K^r_{i(\sigma'_{j+1})}.\]
Let $i\in N^r$ such that $K^r_i\subset  K_{i(\sigma_{j+1})}\setminus K^r_{i(\sigma'_{j+1})}$ and $w^r_{\sigma'_j}\left(K^r_{i(\sigma'_j)}\right)\subset K^r_i$. Since $\M$ is positive, there exists $e\in E^r$ such that $i(e) = i$ and $r(e)=\sigma_{j+1}$. Set $\sigma'_{j+1}:=e$. Then $t(\sigma'_j)=i(\sigma'_{j+1})$ and $\Psi_r(\sigma')=\sigma$. Iterate the procedure until $j$ exceeds the maximal absolute value encountered so far and set $\sigma^1:=\sigma'$. By iterating the procedure, we obtain a sequence $(\sigma^n)_{n\in\N}\subset\bar\Sigma^r$ such that there exists an increasing sequence  $(m_n)_{n\in\N}\subset\N\cup\{0\}$ such that $(\sigma^n_{-m_n},...,\sigma^n_{m_n})$ is a path and $\Psi_r(\sigma^n)=\sigma$ for all $n\in\N$. By the compactness and the metrizability of $\bar\Sigma^r$, there exists a subsequence $(\sigma^{n_k})_{k\in\N}$ and $\sigma''\in\bar\Sigma^r$ such that $\sigma^{n_k}\to\sigma''$ as $k\to\infty$. Hence, since $\Psi_r$ is continuous,  $\Psi_r(\sigma'')=\sigma$. In particular, by the definition of $\Psi_r$, $\sigma''_i\in E^r$ for all $i\in\Z$. Suppose, there exists $j\in\Z$ such that $t(\sigma''_j)\ne i(\sigma''_{j+1})$. Then, by the openness of $_j[\sigma''_j, \sigma''_{j+1}]$, it contains infinitely many of $\sigma^n$, but this contradicts to their construction. Thus $\sigma''\in\Sigma^r_G$. Together with (ii), this completes the proof of (iii).

(iv)
Let $x\in K$ and $ _1[e_1,...,e_n]^+\subset\bar\Sigma^+$. Then, by the definition of $\psi_r$, 
\[{\psi_r}^{-1}\left(\L _1[e_1,...,e_n]^+\right)=\bigcup\limits_{e'_1,...,e'_n\in\bar{E^r},r(e'_i)=e_i}\ _1[e'_1,...,e'_n]^+.\]
Therefore,
 \begin{eqnarray}\label{pmre}
     \psi_r(P^r_x)\left(\L_1[e_1,...,e_n]^+\right)=\sum\limits_{e'_1,...,e'_n\in E^r,r(e'_i)=e_i} p^r_{e'_1}(x)...p^r_{e'_n}\circ
w^r_{e'_{n-1}}\circ...\circ w^r_{e'_1}(x).
 \end{eqnarray}
Now, observe that, by the definition of $p^r_{e'}$'s, there exists at most one $e'_1\in E^r$ with $r(e'_1) = e_1$ and $p^r_{e'_1}(x)>0$. For this $e'_1$, $p^r_{e'_1}(x) = p_{e_1}(x)$ and $w^r_{e'_1}(x) = w_{e_1}(x)$. If there are no such $e'_1$, then the right hand side of \eqref{pmre} is zero and $p_{e_1}(x)=0$ also, since $\sum_{e',r(e')=e} p'_{e'}=p_{e}$ for all $e\in E$. Hence
\[\psi_r(P^r_x)\left(\L_1[e_1,...,e_n]^+\right)=\sum\limits_{e'_2,...,e'_n\in E^r,r(e'_i)=e_i} p_{e_1}(x)...p^r_{e'_n}\circ
w^r_{e'_{n-1}}\circ...\circ w^r_{e'_2}\circ w_{e_1}(x).\] 
Now, by repeating the argument for $x_1:=w_{e_1}(x)$,  $x_2:=w_{e_2}\circ w_{e_1}(x)$, ..., $x_{n-1}:=w_{e_{n-1}}\circ...\circ w_{e_1}(x)$, we obtain
\[\psi_r(P^r_x)\left(\ _1[e_1,...,e_n]^+\right)= p_{e_1}(x)...p_{e_n}\circ
w_{e_{n-1}}\circ...\circ w_{e_1}(x)=P_x\left(_1[e_1,...,e_n]^+\right).\] 
Thus, the claim follows, since the class of the cylinder sets generates the $\sigma$-algebra, is $\cap$-stable and has a countable subset covering $\bar\Sigma^{r+}$.

(v) Let $\mu\in P(\M)$ and $ _{-n}[e_{-n},...,e_n]\subset\bar\Sigma$. Then,  the same way as above,
\begin{eqnarray*}
    &&\Psi_r(\Phi_r(\mu))\left(\L_{-n}[e_{-n},...,e_n]\right)\\
&=&\sum\limits_{e'_{-n},...,e'_n\in E^r,r(e'_i)=e_i} \int p^r_{e'_{-n}}(x)...p^r_{e'_n}\circ
w^r_{e'_{n-1}}\circ...\circ w^r_{e'_{-n}}(x)d\mu(x)\\
&=&\Phi(\mu)\left(\L_{-n}[e_{-n},...,e_n]\right).
 \end{eqnarray*}
Thus, the assertion follows.
\hfill$\Box$

\subsection{Equilibrium states}

Now, we are going to define the main objects on the code space which are useful not only for a description of the invariant measures, but also, combined with another object which will be introduced in subsection \ref{ndss},  allow to control the asymptotic behaviour of the system, most importantly at the boundaries of the atoms of the Markov partition, where the continuity of the system may not be available.

Fix $x_i\in K_i$ for all $i\in N$, and set
\[D:=\left\{\sigma\in\Sigma_G\left|\ \lim\limits_{m\to-\infty}w_{\sigma_0}\circ ... \circ w_{\sigma_m}(x_{i(\sigma_m)})\mbox{ exists}\right.\right\}\]
and
\begin{equation*}
    F(\sigma):=\left\{\begin{array}{cc}
    \lim\limits_{m\to-\infty}w_{\sigma_0}\circ w_{\sigma_{-1}}\circ...\circ w_{\sigma_{m}}(x_{i(\sigma_{m})})&  \mbox{if }\sigma\in D\\
     x_{t(\sigma_0)}& \mbox{ otherwise, }
     \end{array}\right.
     \end{equation*}
for all $\sigma\in\bar\Sigma$. $F:\bar\Sigma\longrightarrow K$ is called the {\it coding map} of the Markov system. Clearly, it is  $\mathcal{F}$-Borel-measurable.  Furthermore, let  $F: P_S(\bar\Sigma)\longrightarrow P(K)$ be given by $ F(\Lambda)(B):=\Lambda(F^{-1}(B))$ for all Borel $B\subset K$ and $\Lambda\in P_S(\bar\Sigma)$.

 Next, set
 \[E(\mathcal{M}):=\left\{\Lambda\in P_S(\bar\Sigma)|\  \Lambda(D) = 1\mbox{ and }E_\Lambda(1_{_1[e]}|\mathcal{F})=p_e\circ F\  \Lambda\mbox{-a.e. for all }e\in E\right\}.\]
We call the members of $E(\M)$ the {\it equilibrium states} of $\M$. 

It will be shown in subsection \ref{esss} that the definition of $E(\mathcal{M})$ naturally extends the notion of {\it equilibrium states} in the thermodynamic sense. Also, an anonymous reviewer pointed out that the condition for members in $E(\mathcal{M})$ is related to the 'conformality', as in \cite{DU}.

Now, we show that the property of $E(\mathcal{M})$ is transferable under the refinement in some cases.
\begin{Definition}
     Let $\M^r:=(K^r_{i(e)},w^r_e,p^r_e)_{e\in E^r}$ be a {\it refinement} of $\M$. Define $D^r$ and the coding map $F_r$ associated with $\M^r$ as above by choosing, for every $j\in N^r$, $x^r_j\in K^r_j$ such that $x^r_j = x_i$ if $x_i\in K^r_j$.
\end{Definition}

If all $w_e|_{K_{i(e)}}$'s are contractions, then, obviously, $D$  and $F|_D$ do not depend on the choice of $x_i$'s.
\begin{Lemma}\label{esrl}
Suppose  $D$ and $F|_D$ do not depend on the choice of $x_i$'s, and $\M^r$ is countable. Then the following holds true.\\
 (i)  $\Psi_r(D^r)\subset D$.\\
 (ii) $\Psi_r(D^r)=D$ if $\M$ is positive and $\Psi_r$ is continuous.\\
 (iii) $F_r(\sigma')= F\circ \Psi_r(\sigma')$ for all $\sigma'\in D^r$.\\
 (iv) $\Psi_r(\Lambda^r)\in E(\M)$ for all $\Lambda^r\in E(\M^r)$.
\end{Lemma}
{\it Proof.}    
(i), (iii) Let $\sigma'\in D^r$. Let $\sigma:=\Psi_r(\sigma')$. Then the following limits exist, and, by the hypothesis, 
  \begin{eqnarray}\label{rle}
    F_r(\sigma')&=&\lim\limits_{m\to-\infty}w^r_{\sigma'_0}\circ...\circ w^r_{\sigma'_m}\left(x^r_{i(\sigma'_m)}\right)=\lim\limits_{m\to-\infty}w_{\sigma_0}\circ...\circ w_{\sigma_m}\left(x^r_{i(\sigma'_m)}\right)\nonumber\\
&=&\lim\limits_{m\to-\infty}w_{\sigma_0}\circ...\circ w_{\sigma_m}\left(x_{i(\sigma_m)}\right)=F\circ \Psi_r(\sigma').
   \end{eqnarray}
This shows (i) and (iii).

(ii) Let $\sigma\in D$. By Lemma \ref{prm} (iii), $\sigma\in\Psi_r(\Sigma^r_G)$, i.e. there exists $\sigma'\in\Sigma^r_G$ such that $\Psi_r(\sigma')=\sigma$. Then 
 \begin{eqnarray*}
    \lim\limits_{m\to-\infty}w_{\sigma_0}\circ...\circ w_{\sigma_m}\left(x_{i(\sigma_m)}\right)&=& \lim\limits_{m\to-\infty}w_{\sigma_0}\circ...\circ w_{\sigma_m}\left(x^r_{i(\sigma'_m)}\right)\\
&=&\lim\limits_{m\to-\infty}w^r_{\sigma'_0}\circ...\circ w^r_{\sigma'_m}\left(x^r_{i(\sigma'_m)}\right).
   \end{eqnarray*}
Hence, $\sigma'\in D^r$, and therefore, $\sigma\in \Psi_r(D^r)$. Together with (i), this shows (ii).

(iv)  Let $\Lambda^r\in E(\M^r)$. Then, by (i), 
\[\Psi_r(\Lambda^r)(D)=\Lambda^r\left(\Psi_r^{-1}(D)\right)\geq \Lambda^r\left(\Psi_r^{-1}(\Psi_r(D^r))\right)\geq\Lambda^r\left(D^r\right)=1.\]
Now, let  $A:=\L_m[e_m,...,e_0]\in\mathcal{F}$. Then, obviously, $\Psi_r^{-1}(A)\in\mathcal{F}^r$, and therefore, by (iii), for $e\in E$,
  \begin{eqnarray*}
      &&\int\limits_{A}1_{_1[e]}d\Psi_r\left(\Lambda^r\right)= \int\limits_{\Psi_r^{-1}(A)}1_{\Psi_r^{-1}(\L_1[e])}d\Lambda^r=\sum\limits_{e'\in E^r,\ r(e')=e}\int\limits_{\Psi_r^{-1}(A)}1_{_1[e']}d\Lambda^r\\
&=&\sum\limits_{e'\in E^r,\ r(e')=e}\int\limits_{\Psi_r^{-1}(A)}p_{e'}^r\circ F_rd\Lambda^r=\int\limits_{\Psi_r^{-1}(A)}p_{e}\circ F\circ \Psi_rd\Lambda^r\\
&=&\int\limits_{A}p_{e}\circ Fd\Psi_r\left(\Lambda^r\right).
   \end{eqnarray*}
Since the set of all such $A$ generates $\mathcal{F}$, is $\cap$-stable and covers $\bar\Sigma$, we conclude that $\Psi_r(\Lambda^r)\in E(\M)$.
\hfill$\Box$

\subsubsection{Thermodynamic equilibrium states}\label{esss}

 Now,  we are going to show that the members of $E(\mathcal{M})$ with finite entropy  which can be computed according to Kolmogorov-Sinai Theorem are exactly the equilibrium states in the thermodynamic sense, which minimise the {\it free energy} of the system, for the following energy function. Set
      \begin{equation*}
    u(\sigma):=\left\{\begin{array}{cc}
    \log  p_{\sigma_1}\circ F(\sigma)&  \mbox{if }\sigma\in D\\
    -\infty& \mbox{ otherwise }
     \end{array}\right.
     \end{equation*}
      for all $\sigma\in\bar\Sigma$ with the definition $\log(0) := -\infty$.  $u$ is called the {\it  energy function} of the Markov system.

\begin{Definition}
For $\Lambda\in P_S(\bar\Sigma)$, set
 \begin{equation}\label{KSf}
    h_S(\Lambda):= H_\Lambda\left((\ _1[e])_{e\in \bar E}|\mathcal{F}\right).
   \end{equation}
 Recall that, by Kolmogorov-Sinai Theorem, $h_S(\Lambda)$ is {\it Shannon-Kolmogorov-Sinai entropy}  if $-\sum_{e\in\bar E} \Lambda\left(\O_1[e]\right)\log  \Lambda\left(\O_1[e]\right)<\infty$.
$\Lambda_0\in P_S(\bar\Sigma)$ is said to be an {\it equilibrium state} for $u$ iff $h_S(\Lambda_0)<\infty$ and
         \begin{equation*}
      h_{S}(\Lambda_0) + \int u d\Lambda_0 = \sup\limits_{\Lambda\in P_S(\bar\Sigma)\mbox{, }h_S(\Lambda)<\infty}\left\{ h_S(\Lambda) + \int u d\Lambda\right\}.
   \end{equation*}
 Let $E(u)\subset  P_S(\bar\Sigma)$ denote the set of all equilibrium states for $u$. 
\end{Definition}

\begin{Lemma}\label{fel}
   Let $\Lambda\in P_S(\bar\Sigma)$ such that $h_S(\Lambda)<\infty$. Then  
 \begin{equation}\label{feie}
    h_S(\Lambda) +\int u d\Lambda \leq 0,
   \end{equation}
and the equality holds if and only if $\Lambda\in E(\mathcal{M})$.
\end{Lemma}
{\it Proof.}    The proof  is an adaptation of  Ledrappier's proof \cite{Le} and that of Lemma 5 in \cite{Wer6}.
Let us abbreviate 
\[g_e:= E_{\Lambda}\left(1_{_1[e]}|\mathcal{F}\right)\]
 for all $e\in\bar E$. If $\int u d\Lambda = -\infty$, then $ h_S(\Lambda) +\int u d\Lambda = -\infty<0$  and \eqref{feie} holds true. Otherwise,
$\Lambda(D)=1$, and therefore, since $\int_{\{g_e=0\}}\ 1_{_1[e]}\log p_e\circ F d\Lambda=0$ for all $e\in E$,
 \begin{eqnarray}\label{lei}
h_S(\Lambda) +\int u d\Lambda& =&-\sum\limits_{e\in\bar E}\int g_e\log g_e d\Lambda+\sum\limits_{e\in\bar E}\int1_{_1[e]}\log p_e\circ F d\Lambda\nonumber\\
& =&\sum\limits_{e\in\bar E}\int g_e\log\frac{p_e\circ F}{g_e} d\Lambda\nonumber\\
& \leq&\sum\limits_{e\in\bar E}\int g_e\left(\frac{p_e\circ F}{g_e}-1\right) d\Lambda\nonumber\\
& =&\sum\limits_{e\in\bar E}\int \left(p_e\circ F- g_e\right) d\Lambda\nonumber\\
&=&0.
\end{eqnarray}
Thus \eqref{feie} holds true in this case also. If \eqref{feie} is an equality, then $\int u d\Lambda > -\infty$, and therefore, $\Lambda(D) = 1$ and, by \eqref{lei},
\[\sum\limits_{e\in\bar E}\int g_e\log\frac{p_e\circ F}{g_e} d\Lambda=\sum\limits_{e\in\bar E}\int g_e\left(\frac{p_e\circ F}{g_e}-1\right) d\Lambda.\]
Hence
\[\log\frac{p_e\circ F(\sigma)}{g_e(\sigma)} = \left(\frac{p_e\circ F(\sigma)}{g_e(\sigma)}-1\right)\mbox{ for }\Lambda\mbox{-a.e. }\sigma\in\left\{g_e>0\right\}\]
for all $e\in\bar E$, but this is possible if and only if
\[ g_e(\sigma) =p_e\circ F(\sigma)\mbox{ for }\Lambda\mbox{-a.e. }\sigma\in\left\{g_e>0\right\}\]
for all $e\in\bar E$. Therefore,
\[  E_{\Lambda}\left(1_{_1[e]}|\mathcal{F}\right) \leq p_e\circ F\ \ \ \Lambda\mbox{-a.e.}\]
for all $e\in\bar E$. However, as then
\[1=\sum\limits_{e\in\bar E}\int E_{\Lambda}\left(1_{_1[e]}|\mathcal{F}\right) d\Lambda\leq \sum\limits_{e\in\bar E}\int p_e\circ F d\Lambda = 1,\]
it follows that
\[  E_{\Lambda}\left(1_{_1[e]}|\mathcal{F}\right) = p_e\circ F\ \ \ \Lambda\mbox{-a.e.}\]
for all $e\in\bar E$. Thus $\Lambda\in E(\mathcal{M})$.

Conversely, if $\Lambda\in E(\mathcal{M})$, then, as $\Lambda(\O _1[\infty]) = 0$,
 \begin{eqnarray*}
 h_S(\Lambda)&=&-\sum\limits_{e\in E}\int E_{\Lambda}\left(1_{_1[e]}|\mathcal{F}\right)\log  E_{\Lambda}\left(1_{_1[e]}|\mathcal{F}\right)d\Lambda\\
&=&-\sum\limits_{e\in E}\sum\limits_{n\leq 0}\int\limits_{\left\{n-1< \log  E_{\Lambda}\left(1_{_1[e]}|\mathcal{F}\right)\leq n\right\}} E_{\Lambda}\left(1_{_1[e]}|\mathcal{F}\right)\log  E_{\Lambda}\left(1_{_1[e]}|\mathcal{F}\right)d\Lambda\\
&=&-\sum\limits_{e\in E}\sum\limits_{n\leq 0}\int\limits_{\left\{n-1< \log  E_{\Lambda}\left(1_{_1[e]}|\mathcal{F}\right)\leq n\right\}} 1_{_1[e]}\log  E_{\Lambda}\left(1_{_1[e]}|\mathcal{F}\right)d\Lambda\\
&=&-\sum\limits_{e\in E}\int\limits_{_1[e]}\log  p_e\circ Fd\Lambda\\
&=&-\int u d\Lambda.
\end{eqnarray*}
That is
\[ h_S(\Lambda) +\int u d\Lambda = 0.\]
This completes the proof.
\hfill$\Box$

\begin{Theorem}\label{est}
    If $\{M\in E(\mathcal{M})|\ h_S(M)<\infty\}$ is not empty, then $\{M\in E(\mathcal{M})|\  h_S(M)<\infty\}= E(u)$.
\end{Theorem}
{\it Proof.}  
By Lemma \ref{fel}, every member of $\{M\in E(\mathcal{M})|\  h_S(M)<\infty\}$ is an equilibrium state of $u$. Conversely,  for every $\Lambda_0\in E(u)$, by the hypothesis and Lemma \ref{fel},
  \[h_{S}(\Lambda_0) + \int u d\Lambda_0 = 0.\]
Thus, by Lemma \ref{fel}, $\Lambda_0\in E(\mathcal{M})$. This completes the proof.
\hfill$\Box$

Theorem \ref{est} and Example \ref{ieex} from Section \ref{easec} seem to indicate that Shannon-Kolmogorov-Sinai entropy might be not the best choice of the entropy for a satisfactory thermodynamic description of such systems.

\subsection{Uniformly continuous Markov system}
In this subsection, we develop a general theory on the relation of the equilibrium states and the invariant Borel probability measures measures of $\mathcal{M}$ if it is uniformly continuous.

 Let $e\in E$. If $w_e|_{K_{i(e)}}$ is uniformly continuous, let $\bar w_e$ denote the continuous extension of $w_e|_{K_{i(e)}}$ on $\bar K_{i(e)}$, which then can be considered to be extended on $K$ arbitrarily.

\begin{Proposition}\label{cml}
    Suppose $w_e|_{K_{i(e)}}$ is uniformly continuous for all $e\in E$.  
Then $F(M)\in P(\mathcal{M})$ for all $M\in E(\mathcal{M})$.
\end{Proposition}
{\it Proof.}  
   Let $M\in E(\mathcal{M})$, and $f\in\mathcal{L}^B(K)$ be bounded. Observe that, since $M(D) = 1$,
 \begin{equation}\label{cmcr}
     \bar w_{\sigma_1}\circ F(\sigma) = F\circ S(\sigma) \mbox{ for }M\mbox{-a.a. }\sigma\in\bar\Sigma.
\end{equation}
 Then, by the shift invariance of $M$,
\begin{eqnarray*}
          &&\int f\circ F\ dM=\int f\circ F\circ S\ dM=\sum\limits_{e\in E}\int 1_{_1[e]}f\circ\bar w_e\circ F\ dM\\
&=&\sum\limits_{e\in E}\int p_e\circ Ff\circ w_e\circ F\ dM=\int Uf\circ F\ dM.
      \end{eqnarray*}
Thus
\[ \int f\ dF(M)=\int f\ dU^{*}F(M)\]
for all bounded $f\in\mathcal{L}^B(K)$.
\hfill$\Box$

Now, for $\mu\in P(\mathcal{M})$, set  
\[\phi_m(\mu)(A):= \int P^m_x(A)d\mu(x)\]
for all $A\in\mathcal{A}_m$ and $m\leq 0$. Observe that, by the invariance of $\mu$,  $\phi_m(\mu)$'s are consistent for all $m\leq 0$ (e.g. see \cite{Wer10}). Let $\Phi(\mu)\in P_S(\bar\Sigma)$ denote the measure which uniquely extends $\phi_m(\mu)$'s on the Borel $\sigma$-algebra, e.g. by the Kolmogorov Consistency Theorem. This defines a map $\Phi: P(\mathcal{M})\longrightarrow P_S(\bar\Sigma)$.

 It is not difficult to check that, for every $\nu\in P(K)$ and $\Omega\in\mathcal{B}(K)\otimes\mathcal{B}(\bar\Sigma^+)$,
\[\tilde\phi(\nu)(\Omega):=\int P_x\left(\{\sigma\in\bar\Sigma^+|\ (x,\sigma)\in\Omega\}\right)d\nu(x)\]
defines a probability measure on product $\sigma$-algebra $\mathcal{B}(K)\otimes\mathcal{B}(\bar\Sigma^+)$ such that
\[\int fd\tilde\phi(\nu)=\int\int f(x,\sigma)dP_x(\sigma)d\nu(x)\]
for every $\tilde\phi(\nu)$-integrable function $f:K\times\bar\Sigma^+\longrightarrow[-\infty,+\infty]$.

Set
\begin{eqnarray*}
\pi:&\bar\Sigma&\longrightarrow\bar\Sigma^+\\
&(...,\sigma_{-1},\sigma_0,\sigma_1,...)&\longmapsto(\sigma_1,\sigma_2...)
\end{eqnarray*}
and
 \begin{eqnarray*}
       \eta:&\bar\Sigma&\longrightarrow K\times\bar\Sigma^+\\
                 &\sigma&\longmapsto (F(\sigma),\pi(\sigma)).
 \end{eqnarray*}

\begin{Lemma}\label{pesl}
Suppose $w_e|_{K_{i(e)}}$ is uniformly continuous for all $e\in E$. Let $M\in E(\mathcal{M})$.\\
    (i) $\Phi(F(M)) = M$.\\
    (ii) $ \eta(M) = \tilde\phi(F(M))$.\\
    (iii) Let $f_e:K\longrightarrow [-\infty,+\infty]$ be Borel  measurable for all $e\in E$ such that  $\sum_{e\in E}1_{_1[e]}|f_e|\circ F\in\mathcal{L}^1(M)$. Then there exists $g\in \mathcal{L}^1(\tilde\phi(F(M)))$ such that
$\int gd\tilde\phi(F(M))=\sum_{e\in E}\int p_ef_edF(M)$
 and 
\[\frac{1}{n}\sum\limits_{k=1}^nf_{\sigma_{k+1}}\circ w_{\sigma_k}\circ ... \circ w_{\sigma_1}(x)\to g(x,\sigma)\]
for $P_x$-a.e. $\sigma\in\bar\Sigma^+$ and for $F(M)$-a.e. $x\in K$, and  in $\mathcal{L}^1(\tilde\phi(F(M)))$.\\
   (iv) For $F(M)$-a.e $x_0\in K$, the sequence of probability measures $(\alpha^{x_0}_n)_{n\in\mathbb{N}}$ on $N$ given by 
                 $\alpha^{x_0}_n(\{j\}):=1/n\sum_{k=1}^n{U^*}^k\delta_{x_0}\left(K_j\right)$ for all $j\in N$ and $n\in\mathbb{N}$
converges in total variation.
\end{Lemma}
{\it Proof.}  
   (i) Let $_1[e_1,...,e_n]\subset\bar\Sigma$ with $(e_1,...,e_n)\in E$. One easily checks  that, by the shift-invariance of $M$ and \eqref{cmcr},
\begin{eqnarray}\label{escp}
    E_M\left(1_{_1[e_1,...,e_n]}|\mathcal{F}\right) (\sigma)=   P^1_{F(\sigma)}\left( _1[e_1,...,e_n]\right)\ \ \ \ \mbox{ for }M\mbox{-a.e. }\sigma\in\bar\Sigma. 
\end{eqnarray}
 Therefore,
\begin{eqnarray*}
    \Phi(F(M)) \left( _1[e_1,...,e_n]\right) &=& \int P^1_x \left( _1[e_1,...,e_n]\right)dF(M)\\
      &=&  \int P^1_{F(\sigma)} \left( _1[e_1,...,e_n]\right)dM(\sigma)\\
      &=& M \left(_1[e_1,...,e_n]\right).
\end{eqnarray*}
Thus, by the shift-invariance of the measures, they  agree on  the class of  cylinder sets of the form $_m[e_m,...,e_{|m|}]$, $m\leq 0$, where $e_m,...,e_{|m|}\in E$. If $e_i=\infty$ for some $m\leq i\leq |m|$, then $\Phi(F(M)) ( _m[e_m,...,e_{|m|}]) =0$ by the definition, and $M ( _m[e_m,...,e_{|m|}]) =0$, as $M(\bar\Sigma\setminus D)=0$.
Thus, the measures  agree on the Borel $\sigma$-algebra, as the class of  cylinder sets of the form $_m[e_m,...,e_{|m|}]$, $e_m,...,e_{|m|}\in\bar E$, $m\leq 0$, plus empty set,  generates the product $\sigma$-algebra, is $\cap$-stable and, obviously, $\bigcup_{e\in \bar E}\O _0[e]=\bar\Sigma$.

   (ii) We only need to check that
\begin{equation}\label{mpce}
 \eta(M)\left(A\times\ _{1}[e_{1},...,e_n]^+\right)=\tilde \phi(F(M))\left(A\times\ _{1}[e_{1},...,e_n]^+\right)
\end{equation}
for all cylinder sets $\ _{1}[e_{1},...,e_n]^+\subset\bar\Sigma^+$ and Borel
$A\subset K$. For such sets,
\begin{eqnarray*}
  \eta(M)\left(A\times\ _{1}[e_{1},...,e_n]^+\right)&=&
 M\left(F^{-1}(A)\cap\  _{1}[e_{1},...,e_n]\right)\\
 &=&\int\limits_{F^{-1}(A)}1_{_{1}[e_{1},...,e_n]}dM,
\end{eqnarray*}
where $_{1}[e_{1},...,e_n]\subset\bar\Sigma$ is the pre-image of
$_{1}[e_{1},...,e_n]^+$ under $\pi$. Clearly,  both sides of \eqref{mpce} are zero if $e_i = \infty$ for some $1\leq i\leq n$. Now, let $e_{1},...,e_n\in E$. Then, by \eqref{escp}, 
\begin{eqnarray*}
 \int\limits_{F^{-1}(A)}1_{_{1}[e_{1},...,e_n]}dM&=&\int\limits_{F^{-1}(A)}P_{F(\sigma)}\left(
_{1}[e_{1},...,e_n]^+\right)dM(\sigma)\\
&=&\int\limits_{A}P_{x}\left(
_{1}[e_{1},...,e_n]^+\right)dF(M)(x),
\end{eqnarray*}
as desired.

(iii)  Set $f_\infty:=0$ and $v(\sigma):= f_{\sigma_1}(F(\sigma))$ for all $\sigma\in\bar\Sigma$. Then,
 \begin{eqnarray*}
\int |v| dM&=&\int\sum\limits_{e\in E} 1_{_1[e]}|f_e|\circ F dM<\infty.
 \end{eqnarray*}
Hence, $v\in\mathcal{L}^1(M)$. Let $\mathcal{I}$ be the $\sigma$-algebra of all shift-invariant Borel subsets of $\bar\Sigma$. Set $\bar v:=E_{M}(v|\mathcal{I})$. Then, by Birkhoff's Ergodic Theorem,
\[\frac{1}{n}\sum\limits_{k=1}^nv\circ S^k\to\bar v\ \ \ M\mbox{-a.e. and in }\mathcal{L}^1(M).\]
Since $M(D)=1$ and $F\circ S^k(\sigma)=\bar w_{\sigma_k}\circ ... \circ \bar w_{\sigma_{1}}\circ F(\sigma)$ for all $\sigma\in D$ and $k\in\mathbb{N}$,
\[\frac{1}{n}\sum\limits_{k=1}^n f_{\sigma_{k+1}}\circ \bar w_{\sigma_k}\circ ... \circ \bar w_{\sigma_{1}}\circ F(\sigma)\to\bar v(\sigma)\ \ \ M\mbox{-a.e. }\sigma\in\bar\Sigma\mbox{ and in }\mathcal{L}^1(M).\]
Set 
\[\bar f_n(x,\sigma):=\frac{1}{n}\sum\limits_{k=1}^n f_{\sigma_{k+1}}\circ \bar w_{\sigma_k}\circ ... \circ \bar w_{\sigma_{1}}(x)\]
for all $x\in K$, $\sigma\in\bar\Sigma^+$ and $n\in\mathbb{N}$.
Then
\[\bar f_n\circ\eta(\sigma)\to\bar v(\sigma)\ \ \ M\mbox{-a.e. }\sigma\in\bar\Sigma\mbox{ and in }\mathcal{L}^1(M).\]
Hence, $\bar v$ is $\eta^{-1}(\mathcal{B}(K)\otimes\mathcal{B}(\bar\Sigma^+))$-measurable. Therefore, by the Factorisation Lemma, there exists a $\mathcal{B}(K)\otimes\mathcal{B}(\bar\Sigma^+)$-measurable function $g$ such that
\[\bar v = g\circ \eta.\]
Then, by (ii) and the definition of $\tilde\phi(F(M))$, 
\[\bar f_n(x, \sigma)\to g(x,\sigma)\ \ \ \mbox{ for }F(M)\mbox{-a.e. }x\in K\mbox{ and }P_x\mbox{-a.e. } \sigma\in\bar\Sigma^+\]
 and in $\mathcal{L}^1(\tilde\phi(F(M)))$, and
 \begin{eqnarray*}
\int g d\tilde\phi(F(M)) &=& \int g\circ \eta dM =  \int v dM =\sum\limits_{e\in E}\int p_e\circ F f_e\circ F dM\\
&=&\sum\limits_{e\in E}\int p_e f_e dF(M).
\end{eqnarray*}
This completes the proof of (iii), as $P_x(\{\sigma\in\bar\Sigma^+|\ x\notin K_{i(\sigma_1)}$ or $\exists k\in\mathbb{N}$ s.t. $i(\sigma_{k+1})\neq t(\sigma_{k})\}) = 0$, by Remark \ref{ppr}.

(iv) By Proposition \ref{cml}, $\mu:=F(M)\in P(\mathcal{M})$.  
Let $i\in N$.  Set $E_i:=\{e\in E|\ i(e)=i\}$ and $f_e=1_{K_i}$ for all $e\in E_i$ and $f_e=0$ for all $e\in E\setminus E_i$. Then, by  (iii), there exists $g_i\in\mathcal{L}^1(\tilde\phi(\mu))$ such that
\begin{equation}\label{lmm}
\int g_i d\tilde\phi(\mu) = \sum_{e\in E}\int p_ef_e d\mu=\mu(K_i)
\end{equation} 
 and
\[\frac{1}{n}\sum\limits_{k=1}^n f_{\sigma_{k+1}}\circ  w_{\sigma_k}\circ ... \circ  w_{\sigma_{1}}(x)\to g_i(x,\sigma)\ \ \  \mbox{ for }P_x\mbox{-a.e. } \sigma\in\bar\Sigma^+\mbox{ and }\mu\mbox{-a.e. }x\in K.\]
One readily checks that $\int f_{\sigma_{k+1}}\circ  w_{\sigma_k}\circ ... \circ  w_{\sigma_{1}}(x) dP_x(\sigma) = U^k(1_{K_i})(x)$ for all $x\in K$. Hence, by Lebesgue's Dominated Convergence Theorem,
\[\frac{1}{n}\sum\limits_{k=1}^nU^k(1_{K_i})(x)\to \int g_i(x,\sigma)dP_x(\sigma)\ \ \ \mbox{ for }\mu\mbox{-a.e. }x\in K\]
for all $i\in N$. As $N$ is countable, also for $\mu$-a.e. $x\in K$,
\[\frac{1}{n}\sum\limits_{k=1}^nU^k(1_{K_i})(x)\to \int g_i(x,\sigma)dP_x(\sigma)\ \ \ \mbox{ for all }i\in N.\]
Hence, by the Fatou Lemma, 
\begin{eqnarray*}
\sum\limits_{i\in N}\int g_i(x,\sigma)dP_x(\sigma)&=&\sum\limits_{i\in N}\liminf\limits_{n\to\infty}\frac{1}{n}\sum\limits_{k=1}^nU^k(1_{K_i})(x)\\
&\leq&\liminf\limits_{n\to\infty}\sum\limits_{i\in N}\frac{1}{n}\sum_{k=1}^nU^k(1_{K_i})(x)\\
&=&1 \mbox{ for } \mu\mbox{-a.e. }x\in K.
\end{eqnarray*} 
Since, by \eqref{lmm}, $\int\sum_{i\in N} \int g_i(x,\sigma)dP_x(\sigma)d\mu(x) = \sum_{i\in N}\mu(K_i) = 1$, there exists a Borel $H\subset K$ with $\mu(H)=1$ such that for every $x_0\in H$,
\[\lim\limits_{n\to\infty}\alpha^{x_0}_n(\{j\})=\alpha^{x_0}(\{j\}) \ \ \ \mbox{ for all }j\in N\]
where  $\alpha^{x_0}$ is the probability measure on $N$ given by $\alpha^{x_0}(\{j\}):= \int g_j(x_0,\sigma)dP_{x_0}(\sigma)$ for all $j\in N$. 
Let  $x_0\in H$. Choose a finite $V_\epsilon\subset N$ such that  $\alpha^{x_0}( V_\epsilon)>1-\epsilon/3$. Then there exist $S\in\mathbb{N}$ such that $\alpha^{x_0}_n( V_\epsilon)>1-\epsilon/3$  for all $n\geq S$ and finite $V_i\subset N$ such that $\alpha^{x_0}_i( V_i)>1-\epsilon/3$ for all $1\leq i<S$. Set $V'_\epsilon:=V_\epsilon\cup\bigcup_{i=1}^{S-1}V_i$. Then $\alpha^{x_0}_n(V'_\epsilon)>1-\epsilon/3$ for all $n\in\mathbb{N}$. Now, choose $n_0\in\mathbb{N}$ such that $\sum_{j\in V'_\epsilon}|\alpha^{x_0}_n(\{j\})-\alpha^{x_0}(\{j\})|<\epsilon/4$ for all $n\geq n_0$. Let $A\subset N$. Then
\begin{eqnarray*}
\left|\alpha^{x_0}_n(A)-\alpha^{x_0}(A)\right|&\leq& \left|\alpha^{x_0}_n(A\cap V'_\epsilon)-\alpha^{x_0}(A\cap V'_\epsilon)\right|+\alpha^{x_0}_n(N\setminus V'_\epsilon)+\alpha^{x_0}(N\setminus V'_\epsilon)\\
&<&\sum\limits_{j\in V'_\epsilon}\left|\alpha^{x_0}_n(\{j\})-\alpha^{x_0}(\{j\})\right|+\frac{2\epsilon}{3}\\
&<&\frac{\epsilon}{4}+\frac{2\epsilon}{3}
\end{eqnarray*}
for all $n\geq n_0$. Thus $\sup_{A\subset N}|\alpha^{x_0}_n(A)-\alpha^{x_0}(A)|<\epsilon$ for all $n\geq n_0$.
\hfill$\Box$

\subsubsection{The non-degeneracy condition}\label{ndss}

Now, we are going to specify the case where the association of the directed graph with the topological structure of $\M$ does not degenerate asymptotically almost surely, i.e.  $F(\sigma)\in K_{t(\sigma_0)}$ for almost every $\sigma\in\Sigma_G$ with respect to every asymptotic state, which now will be defined also.

\begin{Definition}
     Set $T_j:=\{\sigma\in\Sigma_G|\ t(\sigma_0)=j\}$ for all $j\in N$. Then, obviously, $T_j\cap T_{j'} =\emptyset$ for all $j\neq j'$ and $\bigcup_{j\in N}T_j=\Sigma_G$. Suppose $p_e|_{K_{i(e)}}$ is uniformly continuous for all $e\in E$. For each $e\in E$, let $\bar p_e$ denote the continuous extension of $p_e|_{K_{i(e)}}$ on $\bar K_{i(e)}$ which is extended further on $K$ by zero. Let $\tilde E(\mathcal{M})$ denote the set
\[\left\{\Lambda\in P_S(\bar\Sigma)|\  \Lambda(D) = 1\mbox{ and }E_\Lambda(1_{_1[e]}|\mathcal{F})=\bar p_e\circ F1_{T_{i(e)}}\  \Lambda\mbox{-a.e. for all } e\in E\right\}.\]
We call the members of $\tilde E(\mathcal{M})$ the {\it asymptotic states} of $\M$.

Below it will be shown that every equilibrium state is an asymptotic state and that the converse is true only in the case which we are going to specify now.

We call  $\mathcal{M}$ {\it non-degenerate} if and only if  for every $\Lambda\in \tilde E(\mathcal{M})$ there exists $i\in N$ such that $\Lambda(T_i\cap F^{-1}(K_i))>0$. (For example, every uniformly continuous Markov system with an open partition is non-degenerate, as $T_i\subset F^{-1}(\bar K_i)$ for all $i\in N$.) Set $G:=\bigcup_{k=0}^\infty\bigcup_{i\in N}S^{-k}(F^{-1}(K_i)\cap T_i)$.  Then obviously, $\mathcal{M}$ is non-degenerate if and only if $M(G)>0$ for all  $M\in\tilde E(\mathcal{M})$. For $M\in\tilde E(\mathcal{M})$, set 
 \begin{equation*}
    M_G(B):=\left\{\begin{array}{cc}
    \frac{M(B\cap G)}{M(G)}&  \mbox{ if }M(G)>0\\
    0& \mbox{otherwise }
     \end{array}\right.\mbox{ for all }B\in \mathcal{B}(\bar\Sigma).
     \end{equation*}
Then, clearly, $M_G\in P_S(\bar\Sigma)\cup\{0\}$, as $S^{-1}G\subset G$. Set
$\partial p_e:=\bar p_{e}1_{\bar K_{i(e)}\setminus K_{i(e)}}$ for all $e\in E$, and let $ E_\perp(\mathcal{M})$ denote the set
\[\left\{\Lambda\in P_S\left(\bar\Sigma\right)\left|\ \Lambda(D)=1\mbox{ and } E_{\Lambda}\left(1_{_1[e]}|\mathcal{F}\right) =  \partial p_e\circ F1_{T_{i(e)}}\ \Lambda\mbox{-a.e. for all }\right. e\in E\right\}.\]
\end{Definition}

\begin{Lemma}\label{esfg}
(i)   Let $\Lambda\in E(\mathcal{M})$. Then $\Lambda(G)=1$.\\
(ii)  Suppose $p_e|_{K_{i(e)}}$ is uniformly continuous for all $e\in E$. Then $\Lambda(G)=0$ for all $\Lambda\in E_\perp(\mathcal{M})$. 
\end{Lemma}
{\it Proof.} 
  (i) Let $i\in N$ and $A\in\mathcal{F}$. Then, since $\Lambda(\Sigma_G)=1$,
\[\Lambda\left(A\cap T_i\right)=\sum\limits_{e\in E, i(e) =i}\int\limits_A1_{_1[e]}d\Lambda=\sum\limits_{e\in E, i(e) =i}\int\limits_Ap_e\circ Fd\Lambda = \Lambda\left(A\cap F^{-1}(K_i)\right).\]
Hence, $\Lambda(T_i) =\Lambda(F^{-1}(K_i))$ and $\Lambda(F^{-1}(K_i)\cap T_i) = \Lambda(F^{-1}(K_i))$ for all $i\in N$.
Therefore,
\[\Lambda(G)\geq \Lambda\left(\bigcup\limits_{i\in N}F^{-1}(K_i)\cap T_i\right)= \sum\limits_{i\in N}\Lambda\left(F^{-1}(K_i)\right)=1.\]

(ii) Let $\Lambda\in E_\perp(\mathcal{M})$, $i\in N$ and $A\in\mathcal{F}$. First, observe that, by the Fatou Lemma, 
\[\sum\limits_{e\in E,i(e)=i}\bar p_e\leq 1_{\bar K_i}.\] 
Therefore, since $\Lambda(\Sigma_G)=1$ and $T_i\subset F^{-1}(\bar K_i)$,
\[\Lambda\left(A\cap T_i\right)=\sum\limits_{e\in E, i(e) =i}\int\limits_A\bar p_e\circ F1_{\bar K_i\setminus K_i}\circ F1_{T_i}d\Lambda \leq\int\limits_A1_{T_i\setminus F^{-1}(K_i)}d\Lambda.\]
Hence $\Lambda( F^{-1}(K_i)\cap T_i)=0$, and therefore, 
\[\Lambda(G)\leq\sum\limits_{i\in N,k\geq 0}\Lambda\left(S^{-k}\left( F^{-1}(K_i)\cap T_i\right)\right)=0.\]
\hfill$\Box$

\begin{Lemma}\label{endl}
      Suppose $\mathcal{M}$ is uniformly continuous.   Let $M\in\tilde E(\mathcal{M})$.\\
(i) If $M(G)>0$, then $M_G\in E(\mathcal{M})$.\\
(ii) $E(\mathcal{M})\subset\tilde E(\mathcal{M})$.\\
(iii) There exist  $\Lambda\in E(\mathcal{M})\cup\{0\}$ and $\Lambda_\perp\in E_\perp(\mathcal{M})\cup\{0\}$ such that
\[M=M(G)\Lambda+(1-M(G))\Lambda_\perp.\]
The decomposition is unique if $0<M(G)<1$.
(Then,  by Lemma \ref{esfg},  $\Lambda = M_G$, and $\Lambda_\perp$ is singular to $\Lambda$.)\\
(iv) $E_\perp(\mathcal{M})\subset\tilde E(\mathcal{M})$.
\end{Lemma}
{\it Proof.} 
(i) Clearly, $M_G\ll M$. Hence, $M_G(D) = 1$. Let $e\in E$. It is a well-known fact and it can be easily checked that the absolute continuity relation and the shift-invariance of the measures imply (the shift-invariance of the Radon-Nikodym derivative which in turn implies) that  $E_M\left(1_{_1[e]}|\mathcal{F}\right) = E_{M_G}\left(1_{_1[e]}|\mathcal{F}\right)$ $M_G$-a.e. Hence, 
\[E_{M_G}\left(1_{_1[e]}|\mathcal{F}\right)(\sigma)= \bar p_e\circ F(\sigma)1_{K_{i(e)}}\left(x_{t(\sigma_0)}\right)\ \ \ \mbox{ for }M_G\mbox{-a.a. }\sigma\in\bar\Sigma.\] 
Let  $\sigma\in G\cap D$. Then there exist $k\geq 0$ and $j\in N$ such that $S^k(\sigma)\in F^{-1}(K_j)\cap T_j$. That is $F\circ S^k(\sigma) \in K_j$ and  $i(\sigma_{k+1})=t(\sigma_k)=j$. Hence,
 $ F\circ S^n(\sigma) = \bar w_{\sigma_n}\circ ... \circ \bar w_{\sigma_{k+1}}(F\circ S^{k}(\sigma))= w_{\sigma_n}\circ ... \circ w_{\sigma_{k+1}}( F\circ S^{k}(\sigma))$  for all $n\geq k$, and therefore,
\begin{eqnarray*}
  \bar p_e\circ F\circ S^n(\sigma) 1_{T_{i(e)}}\circ S^n(\sigma) &=& \bar p_e\circ w_{\sigma_n}\circ ... \circ   w_{\sigma_{k+1}}\left( F\circ S^{k}(\sigma)\right)1_{K_{i(e)}}\left(x_{t(\sigma_n)}\right)\\
& =& p_e\circ w_{\sigma_n}\circ ... \circ  w_{\sigma_{k+1}}\left( F\circ S^k(\sigma)\right)\\
&=&p_e\circ F\circ S^n(\sigma) 
\end{eqnarray*}
  for all $n\geq k$. Let $A\in\mathcal{F}$. Then
 \begin{equation}\label{ndpe}
     \frac{1}{n}\sum\limits_{i=1}^n \left(\bar p_e\circ F\circ S^i(\sigma) 1_{T_{i(e)}}\circ S^i(\sigma)1_{A}\circ S^i(\sigma)- p_e\circ F\circ S^i(\sigma) 1_{A}\circ S^i(\sigma)\right)\to 0
\end{equation}
for all $\sigma\in G\cap D$. By Birkhoff's Ergodic Theorem, there exist $\bar v_e, v_e\in\mathcal{L}^1(M_G)$ such that
$\int_{A}\bar p_e\circ F 1_{T_{i(e)}}dM_G=\int\bar v_edM_G$, $\int_A p_e\circ F dM_G=\int v_edM_G$,
\[\frac{1}{n}\sum\limits_{i=1}^n \bar p_e\circ F\circ S^i1_{T_{i(e)}}\circ S^i1_{A}\circ S^i\to \bar v_e\mbox{ and }\frac{1}{n}\sum\limits_{i=1}^n p_e\circ F\circ S^i1_{A}\circ S^i\to  v_e\]
both $M_G$-a.e. Hence, since $M_G(G\cap D)=1$, $\bar v_e = v_e$ $M_G$-a.e., and therefore,
\[\int\limits_{A}\bar p_e\circ F 1_{T_{i(e)}}dM_G=\int\limits_{A}p_e\circ F dM_G.\]
Thus
 \begin{equation*}
     E_{M_G}\left(1_{_1[e]}|\mathcal{F}\right) =   p_e\circ F\ \ \ M_G\mbox{-a.e.}
\end{equation*}
This completes the proof of (i).

(ii) Let $\Lambda\in E(\mathcal{M})$. Then, by Lemma \ref{esfg}(i), $\Lambda(G\cap D)=1$. Therefore, by  \eqref{ndpe}, the same way as above,
\[\int\limits_A1_{_1[e]}d\Lambda=\int\limits_Ap\circ Fd\Lambda=\int\limits_A\bar p\circ F1_{T_{i(e)}}d\Lambda\]
for all $e\in E$ and $A\in\mathcal{F}$. Thus $\Lambda\in\tilde E(\mathcal{M})$. This completes the proof of (ii).

(iii) By (i), we can assume that $M(G)<1$. Set \[\Lambda_\perp(B):=\frac{M\left(B\cap \bar\Sigma\setminus G\right)}{M\left(\bar\Sigma\setminus G\right)}\] for all $B\in \mathcal{B}(\bar\Sigma)$. Then $\Lambda_\perp\in P_S(\bar\Sigma)$, since $S^{-1}G\subset G$, and 
 \begin{equation*}
M=M(G)M_G+(1-M(G))\Lambda_\perp.
\end{equation*}
 By (i), $M_G\in E(\mathcal{M})$. Note that $\Lambda_\perp\ll M$. Let $e\in E$. Then, $\Lambda_\perp(D)=1$, and,  as in the proof of (i),
\[E_{\Lambda_\perp}\left(1_{_1[e]}|\mathcal{F}\right)= \bar p_e\circ F1_{T_{i(e)}}\ \ \ \Lambda_\perp\mbox{-a.e. }\] 
 Let $A\in\mathcal{F}$.  By Birkhoff's Ergodic Theorem, there exist $\bar v_e,\partial v_e\in\mathcal{L}^1(\Lambda_\perp)$ such that $\int\bar v_ed\Lambda_\perp =\int_A\bar p_e\circ F1_{T_{i(e)}}d\Lambda_\perp$ and $\int\partial v_ed\Lambda_\perp =\int_A\partial p_e\circ F1_{T_{i(e)}}d\Lambda_\perp$
and
\[\frac{1}{n}\sum\limits_{k=1}^n \bar p_e\circ F\circ S^k1_{T_{i(e)}}\circ S^i1_{A}\circ S^k\to \bar v_e\mbox{ and }\frac{1}{n}\sum\limits_{k=1}^n \partial p_e\circ F\circ S^k1_{A}\circ S^k\to  \partial v_e\]
both $\Lambda_\perp$-a.e.
 Note that $\bar p_e=p_e+\partial p_e$. Therefore,
\[\bar p_e\circ F1_{T_{i(e)}}=p_e\circ F1_{F^{-1}\left(K_{i(e)}\right)\cap T_{i(e)}}+\partial p_e\circ F1_{T_{i(e)}}.\]
Let $\sigma\in D\cap\bar\Sigma\setminus G$. Then for each $k\in\mathbb{N}\cup\{0\}$ and $i\in N$, $S^k(\sigma)\in\bar\Sigma\setminus (F^{-1}(K_i)\cap T_i)$. Hence,
 \begin{equation}\label{sppe}
   \bar p_e\circ F\circ S^k(\sigma)1_{T_{i(e)}}\circ S^k(\sigma)=\partial p_e\circ F\circ S^k(\sigma)1_{T_{i(e)}}\circ S^k(\sigma)
\end{equation}
for all $k\in\mathbb{N}$. Therefore, since $\Lambda_\perp(D\cap\bar\Sigma\setminus G)=1$, $\bar v_e = \partial v_e$  $\Lambda_\perp$-a.e., and
\[\int\limits_A\bar p_e\circ F1_{T_{i(e)}}d\Lambda_\perp=\int\limits_A\partial p_e\circ F1_{T_{i(e)}}d\Lambda_\perp.\]
Thus $\Lambda_\perp\in E_{\perp}(\mathcal{M})$, and the existence of the decomposition is proved.

Suppose $0<M(G)<1$, and there exist  $\Lambda'\in E(\mathcal{M})$ and $\Lambda'_\perp\in E_\perp(\mathcal{M})$ such that
\[M=M(G)\Lambda'+(1-M(G))\Lambda'_\perp.\]
Then, by Lemma \ref{esfg}(i),  $\Lambda'(G)=1$. Hence $\Lambda'_\perp(G) = 0$. Therefore, $M(B\cap G) = M(G)\Lambda'(B)$ for all $B\in\mathcal{B}(\bar\Sigma)$. That is $\Lambda' = M_G$. Hence,  $\Lambda'_\perp = (M-M(G)M_G)/(1-M(G)) = \Lambda_\perp$.
Thus, the decomposition is unique.

(iv) Let $\Lambda\in E_\perp(\mathcal{M})$, $e\in E$ and $A\in\mathcal{F}$. Then, by Lemma \ref{esfg}(ii), $\Lambda(D\cap\bar\Sigma\setminus G)=1$.  Hence, by \eqref{sppe}, the same way as above,
\[\int\limits_A1_{_1[e]}d\Lambda=\int\limits_A\partial p_e\circ F1_{T_{i(e)}}d\Lambda=\int\limits_A\bar p\circ F1_{T_{i(e)}}d\Lambda.\]
Thus $\Lambda\in\tilde E(\mathcal{M})$.
\hfill$\Box$

\begin{Theorem}\label{cndc}
Suppose $\mathcal{M}$ is uniformly continuous. Then the following are equivalent.\\
(i) $\mathcal{M}$ is non-degenerate.\\
(ii) $M(G)>0$ for all $M\in \tilde E(\mathcal{M})$.\\
(iii) $E_\perp(\mathcal{M})$ is empty.\\
(iv) $\tilde E(\mathcal{M})=E(\mathcal{M})$.\\
(v) $M(G)=1$ for all $M\in \tilde E(\mathcal{M})$.
\end{Theorem}
{\it Proof.} 
 $(i)\Leftrightarrow(ii)$ is obvious.

$(ii)\Rightarrow(iii)$ follows by Lemma \ref{endl} (iv) and Lemma \ref{esfg} (ii).

$(iii)\Rightarrow(iv)$  follows by Lemma  \ref{endl} (ii) and (iii).

$(iv)\Rightarrow(v)$  follows by Lemma  \ref{esfg} (i).

$(v)\Rightarrow(ii)$ is obvious.
\hfill$\Box$

\subsubsection{A sufficient condition for the non-degeneracy}
The following lemma will be used to show the non-degeneracy of most of the examples in this article.

\begin{Definition}\label{dod}
Suppose  $\mathcal{M}$ is uniformly continuous. Set 
\[Rf:=\sum\limits_{e\in E}\partial p_e f\circ\bar w_e\]
for all Borel-measurable $f:K\longrightarrow [0,+\infty]$, and
\[\Omega:=\bigcap\limits_{n\in\mathbb{N}}\left\{R^n1\geq1\right\}.\] 
 Clearly, if the partition of $\mathcal{M}$ consists of open sets,  then $R1=0$,  and therefore $\Omega=\emptyset$. Note  that, by the Fatou Lemma,
\[R1=\sum\limits_{e\in E}\partial p_e = \sum\limits_{j\in N}\sum\limits_{e\in E, i(e)=j}\bar p_e1_{\bar K_j\setminus K_j}\leq \sum\limits_{j\in N}1_{\bar K_j\setminus K_j}.\]
Therefore, $\Omega\subset\bigcup_{j\in N}\bar K_j\setminus K_j$.
\end{Definition}

\begin{Lemma}\label{Rpl}
  Suppose $\mathcal{M}$ is uniformly continuous. Let $\Lambda\in E_{\perp}(\mathcal{M})$. Then \[F(\Lambda)(\Omega)=1.\]
\end{Lemma}
{\it Proof.}
 Let  $n\in\mathbb{N}$.
Using the shift-invariance of $\Lambda$ and \eqref{cmcr},  one easily checks that,  for all $e_1,...,e_n\in E$,
\begin{eqnarray*}
&& E_{\Lambda}\left(1_{_1[e_1,...,e_n]}|\mathcal{F}\right)\\
&=& 1_{K_{i(e_2)}}(x_{t(e_1)}) ...1_{K_{i(e_n)}}(x_{t(e_{n-1})})\\
&&\times \partial p_{e_1}\circ F1_{T_{i(e_1)}}\partial p_{e_2}\circ\bar w_{e_1}\circ F...\partial p_{e_n}\circ \bar w_{e_{n-1}}\circ ...\circ\bar w_{e_1}\circ F\ \ \ \ \Lambda\mbox{-a.e.}
\end{eqnarray*}
Therefore, for every $B\in\mathcal{B}(K)$,
\begin{eqnarray*}
F(\Lambda)\left( B\right)&=&\sum\limits_{e_1,...,e_n\in E}\int\limits_{F^{-1}\left(B\right)}  1_{_1[e_1,...,e_n]}d\Lambda
\leq\int\limits_{F^{-1}\left(B\right)}\left(R^{n}1\right)\circ Fd\Lambda\\
&=&\int\limits_{B}R^{n}1dF(\Lambda).
\end{eqnarray*}
That is, $1\leq R^n1$ $F(\Lambda)$-a.e. The assertion follows.
 \hfill$\Box$

\begin{Lemma}\label{ndcl}
  Suppose $\mathcal{M}$ is uniformly continuous. Then $\mathcal{M}$ is non-degenerate if $F^{-1}(\Omega)$ is empty.
\end{Lemma}
{\it Proof.}
 The assertion follows immediately from Lemma \ref{Rpl} and Theorem \ref{cndc}.
 \hfill$\Box$

\subsubsection{The consistency condition}\label{ccs}
 Example  \ref{dMse} (below) shows that the non-degeneracy is not a necessary condition for the existences of an invariant measure for a finite, contractive and uniformly continuous Markov system. Now, we are going to weaken the non-degeneracy condition, in order to include the case where a degenerate system still has an invariant measure because every degeneracy with respect to an asymptotic state transpires in a consistent way, so that no loss of measure occurs. For example, this can happen because of some continuities at some boundaries of some atoms of the partition, or even, as in Example  \ref{dMse}, because the degenerate system is actually a refinement of a non-degenerate one. Of course, in the latter case, one might want to consider instead the non-degenerate Markov system associated with the random dynamical system,  but it might be not favourable because of the loss of some other nice properties, such as the contraction on average. The  theorem below can be used in a degenerate case. 

\begin{Definition}\label{cd}
Suppose  $\mathcal{M}$ is uniformly continuous. We call $\mathcal{M}$ {\it consistent}  if and only if
$F(M)\in P(\mathcal{M})$ for all $M\in \tilde E(\mathcal{M})$. 
\end{Definition}
By Theorem \ref{cndc}  and Proposition \ref{cml}, every uniformly continuous, non-degenerate Markov system is consistent.

\begin{Condition}\label{csc}
 $\mathcal{M}$ is uniformly continuous and
\[\left(Rf\right)\circ F(\sigma) \leq\left(Uf\right)\circ F(\sigma)\mbox{ for all }\sigma\in F^{-1}(\Omega)\]
and all bounded  $f\in\mathcal{L}^B(K)$.
\end{Condition}
Obviously, the condition is satisfied if $F^{-1}(\Omega)$ is empty, which, by Lemma \ref{ndcl}, implies the non-degeneracy. In general, it implies the consistency, as the next theorem shows. However, it is not a necessary condition for it  (see Example  \ref{dMse}). 

\begin{Theorem}\label{imct}
   (i) Uniformly continuous $\mathcal{M}$  is consistent if and only if $F(\Lambda)\in  P(\mathcal{M})$ for all $\Lambda\in  E_\perp(\mathcal{M})$.\\
   (ii) $\mathcal{M}$  is consistent if it satisfies Condition \ref{csc}.
\end{Theorem}
{\it Proof.}
(i) The 'only if' part follows by Lemma \ref{endl} (iv). For the 'if' part, let $M\in \tilde E(\mathcal{M})$ and $\Lambda\in  E_\perp(\mathcal{M})$ such that $M=M(G)M_G+(1-M(G))\Lambda$, by Lemma \ref{endl} (iii). Then $F(M)=M(G)F(M_G)+(1-M(G))F(\Lambda)$. Hence, by Lemma \ref{endl} (i) and Proposition \ref{cml},  $F(\Lambda)\in P(\mathcal{M})$. 

(ii) Let $\Lambda\in  E_\perp(\mathcal{M})$ and $f\in\mathcal{L}^B(K)$ be bounded. Then, by the hypothesis and Lemma \ref{Rpl},
\begin{eqnarray*}
&&\int f dF(\Lambda)=\sum\limits_{e\in E}\int 1_{_0[e]}f\circ Fd\Lambda=\sum\limits_{e\in E}\int 1_{_1[e]}f\circ \bar w_e\circ Fd\Lambda\\
&=&\sum\limits_{e\in E}\int \partial p_e\circ F1_{T_{i(e)}}f\circ \bar w_e\circ Fd\Lambda\leq\int (Rf)\circ F d\Lambda\leq\int Uf dF(\Lambda).
\end{eqnarray*}
Hence
\[\int f dF(\Lambda)\leq \int f dU^*(F(\Lambda)).\]
Since $f$ was arbitrary, it follows that $F(\Lambda)\in P(\mathcal{M})$.
 \hfill$\Box$

\begin{Lemma}\label{uccl}
  Suppose $\mathcal{D}_R$ is  continuous, and $\mathcal{M}$ is uniformly continuous.\\
(i) Let  $f\in\mathcal{L}^B(K)$ be bounded. Then 
 \[Rf\leq\left(\sum\limits_{i\in N}1_{\bar K_i\setminus K_i}\right)Uf.\]
(ii) $\mathcal{M}$ is consistent if $\sum_{i\in N}1_{\bar K_i\setminus K_i}(F(\sigma))\leq 1$ for all $\sigma\in F^{-1}(\Omega)$.
\end{Lemma}
{\it Proof.}
(i)  Let $\mathcal{D}_R=(K, w'_e,p'_e)_{e\in E'}$ with $w'_e:K\longrightarrow K$ and $p'_e:K\longrightarrow[0,1]$ both continuous for all $e\in E'$ and $c:E\longrightarrow E'$ be given  by $w'_{c(e)}|_{K_{i(e)}}=w_e|_{K_{i(e)}}$ and $p'_{c(e)}|_{K_{i(e)}}=p_e|_{K_{i(e)}}$ for all $e\in E$. Note that, for each $j\in N$,
\[\sum\limits_{e_0\in E, i(e_0)=j}\sum\limits_{e\in c^{-1}(\{c(e_0)\})}p_ef\circ w_e\leq Uf.\]
 Furthermore, by the continuity of $\mathcal{D}_R$, for each $e_0\in E$,
\begin{eqnarray*}
\partial p_{e_0}&=&\bar p_{e_0}1_{\bar K_{i(e_0)}\setminus K_{i(e_0)}}=1_{\bar K_{i(e_0)}\setminus K_{i(e_0)}}p'_{c(e_0)}\\
&=&1_{\bar K_{i(e_0)}\setminus K_{i(e_0)}}\sum\limits_{e\in  c^{-1}(\{c(e_0)\})}p_e1_{K_{i(e)}}, 
\end{eqnarray*}
and $\bar w_{e_0}|_{K_{i(e)}\cap\bar K_{i(e_0)} }=w'_{c(e_0)}|_{K_{i(e)}\cap\bar K_{i(e_0)} }=w_e|_{K_{i(e)}\cap\bar K_{i(e_0)} }$ for all
$e\in  c^{-1}(\{c(e_0)\})$.
Therefore,
\begin{eqnarray*}
Rf&=&\sum\limits_{e_0\in E}\partial p_{e_0}f\circ\bar w_{e_0}\\
&=&\sum\limits_{j\in N}\sum\limits_{e_0\in E, i(e_0)=j}1_{\bar K_j\setminus K_j}\sum\limits_{e\in  c^{-1}(\{c(e_0)\})}p_e1_{K_{i(e)}}f\circ\bar w_{e_0}\\
&=&\sum\limits_{j\in N}1_{\bar K_j\setminus K_j}\sum\limits_{e_0\in E, i(e_0)=j}\sum\limits_{e\in  c^{-1}(\{c(e_0)\})}p_e1_{K_{i(e)}}f\circ w_{e}\\
&\leq&\sum\limits_{j\in N}1_{\bar K_j\setminus K_j}Uf.
\end{eqnarray*}
This completes the proof of (i).

(ii) The assertion follows immediately from (i) by Theorem \ref{imct} (ii) and the hypothesis.
\hfill$\Box$

\begin{Remark}
      The consistency condition is also not a necessary condition for the existences of an invariant measure for a finite, contractive and uniformly continuous Markov system. For example, an inconsistent system still can have an invariant measure because there exists $M\in \tilde E(\mathcal{M})$ such that $F(M)\in P(\mathcal{M})$, e.g. see Example \ref{gnce}. However, this situation can be often reduced, as in Example \ref{gnce}, to the consistency of a subsystem, see Corollary \ref{sndct}.
\end{Remark}

\subsubsection{The dominating Markov chain}

The next object  arises naturally from the following condition, which will be used for several purposes in the case of an infinite  set $\{e\in E|\ i(e) = j\}$, $j\in N$.
\begin{Definition}\label{dmcd}
      We say that $\mathcal{M}$ has a {\it dominating} Markov chain iff for every $i\in N$,
\begin{equation}\label{dmcc}
       \xi_i:=\sum\limits_{e\in E, i(e)=i}\sup\limits_{x\in K_i}p_e(x)<\infty.
 \end{equation} 
In this case, set 
\[q_{ij}:=\frac{1}{\xi_i} \sum\limits_{e\in E, i(e)=i, t(e) =j}\sup\limits_{x\in K_i}p_e(x)\]
for all for all $i,j\in N$, and
\begin{equation}\label{usmcc}
       c:=\sum_{j\in N}\sup\limits_{i\in N}\xi_iq_{ij}.
 \end{equation} 
\end{Definition}

\begin{Lemma}\label{acc}
      Suppose $\mathcal{M}$ has a dominating Markov chain and each $p_e|_{K_{i(e)}}$ is uniformly continuous. Let $j\in N$. Then
\[\sum\limits_{e\in E, i(e)=j}\bar p_e(x) = 1_{ \bar K_j}(x)\ \ \ \mbox{ for all }x\in K. \]
\end{Lemma}
{\it Proof.} 
 Let $x\in K$. If $x\notin \bar K_j$ , then, clearly, 
 $\sum_{e\in E, i(e)=j}\bar p_e(x) =0$.
Otherwise, there exists a sequence $(x_n)_{n\in\mathbb{N}}\subset K_j$ such that $\lim_{n\to\infty}x_n = x$. Clearly,   $\sum_{e\in E, i(e)=j}\bar p_e(x_n) =\sum_{e\in E, i(e)=j} p_e(x_n) = 1$ for all $n\in\mathbb{N}$. Since $\lim_{n\to\infty} p_e(x_n) =\bar p_e(x)$ for all $e\in E$ with $i(e)=j$, and $\mathcal{M}$ has a dominating Markov chain, it follows, by Lebesgue's Dominated Convergence Theorem, that 
\[ \sum\limits_{e\in E, i(e)=j}\bar p_e(x) = \lim_{n\to\infty}\sum\limits_{e\in E, i(e)=j} p_e(x_n) =  1.\] 
Thus, combining both cases,
$ \sum_{e\in E, i(e)=j}\bar p_e(x) =1_{\bar K_j}(x) $. 
\hfill$\Box$

\begin{Theorem}\label{ucct}
       Suppose $\mathcal{D}_R$ is  continuous, and $\mathcal{M}$ is uniformly continuous and has a dominating Markov chain. Then $\mathcal{M}$ is consistent.
\end{Theorem}
{\it Proof.} 
 Let $\Lambda\in  E_\perp(\mathcal{M})$ and $f\in\mathcal{L}^B(K)$ be bounded.  Let $\mathcal{D}_R=(K, w'_{e'},p'_{e'})_{e'\in E'}$ where $w'_{e'}:K\longrightarrow K$ and $p'_{e'}:K\longrightarrow[0,1]$ both continuous for all $e'\in E'$ and $c:E\longrightarrow E'$ be given  by $w'_{c(e)}|_{K_{i(e)}}=w_e|_{K_{i(e)}}$ and $p'_{c(e)}|_{K_{i(e)}}=p_e|_{K_{i(e)}}$ for all $e\in E$. Then, by the continuity of $\mathcal{D}_R$, for each $e\in E$, $\partial p_{e}=\bar p_{e}1_{\bar K_{i(e)}\setminus K_{i(e)}}=1_{\bar K_{i(e)}\setminus K_{i(e)}}p'_{c(e)}$ and $\bar w_{e}|_{\bar K_{i(e)} }=w'_{c(e)}|_{\bar K_{i(e)} }$. Note that, by Lemma \ref{acc},  $ \sum_{e\in E, i(e)=j}\partial p_e =1_{\bar K_j\setminus K_j}$ for all $j\in N$. Therefore,  as in the proof of Theorem \ref{imct} (ii),
\begin{eqnarray*}
&&\int f dF(\Lambda)\\
&=&\sum\limits_{e\in E}\int \partial p_e\circ F1_{T_{i(e)}}f\circ \bar w_e\circ Fd\Lambda\\
&=&\sum\limits_{e\in E}\int1_{\bar K_{i(e)}\setminus K_{i(e)}}\circ F p'_{c(e)}\circ F1_{T_{i(e)}}f\circ w'_{c(e)}\circ Fd\Lambda\\
&=&\sum\limits_{e'\in E'}\sum\limits_{e\in E, c(e)=e'}\int p'_{e'}\circ Ff\circ w'_{e'}\circ F1_{\bar K_{i(e)}\setminus K_{i(e)}}\circ F 1_{T_{i(e)}}d\Lambda\\
&=&\sum\limits_{e'\in E'}\sum\limits_{e\in E, c(e)=e'}\sum\limits_{e_0\in E, i(e_0)=i(e)}\int p'_{e'}\circ Ff\circ w'_{e'}\circ F\partial p_{e_0}\circ F 1_{T_{i(e_0)}}d\Lambda\\
&=&\sum\limits_{e'\in E'}\sum\limits_{e\in E, c(e)=e'}\sum\limits_{e_0\in E, i(e_0)=i(e)}\int p'_{e'}\circ Ff\circ w'_{e'}\circ F1_{_1[e_0]}d\Lambda\\
&\leq&\sum\limits_{e'\in E'}\int p'_{e'}\circ Ff\circ w'_{e'}\circ Fd\Lambda\\
&=&\int UfdF(\Lambda).
\end{eqnarray*}
Hence, $\int f dF(\Lambda)\leq \int fdU^*F(\Lambda)$. Since $f$ was arbitrary,  it follows that $U^*F(\Lambda) = F(\Lambda)$. Thus, the assertion follows by Theorem \ref{imct} (i).
\hfill$\Box$

\begin{Lemma}\label{rcl}
Suppose $\mathcal{M}$ is uniformly continuous and positive such that  $D$ and $F|_D$ do not depend on the choice of $x_i$'s. Let $\M^r:=(K^r_{i(e)},w^r_e,p^r_e)_{e\in E^r}$ be a countable refinement of $\M$ which has a dominating Markov chain. Then
\[\Psi_r\left(\Lambda^r\right)\in\tilde E(\M)\ \ \ \mbox{ for all }\Lambda^r\in \tilde E(\M^r).\]
\end{Lemma}
{\it Proof.} By Lemma \ref{endl} and Lemma \ref{esrl} (iv), it is sufficient to show that $\Psi_r(\Lambda^r)\in\tilde E(\M)$ for all $\Lambda^r\in  E_\perp(\M^r)$. So, let $\Lambda^r\in E_\perp(\M^r)$. By Lemma \ref{esrl} (i), $D^r\subset\Psi_r^{-1}(\Psi_r(D^r))\subset\Psi_r^{-1}(D)$. Hence,
\[\Psi_r\left(\Lambda^r\right)(D)=1.\]
 Let $e\in E$. Observe that $\bar p^r_{e'}1_{\bar K^r_{i(e')}\setminus K^r_{i(e')}}=\bar p_{e}1_{\bar K^r_{i(e')}\setminus K^r_{i(e')}}$ for all $e'\in E^r$ with $r(e')=e$, and, since $\M^r$ has a dominating Markov chain, 
\[1_{\bar K^r_{i(e')}\setminus K^r_{i(e')}}=\sum\limits_{e''\in E^r,\ i(e'')=i(e')}\partial p^r_{e''}\]
for all $e'\in E^r$. Let $A:=\L_{n}[e_n,...,e_0]\in\mathcal{F}$. Then, obviously, $\Psi_r^{-1}(A)\in\mathcal{F}^r$, and therefore,
\begin{eqnarray*}
  && \int\limits_A1_{_1[e]}d\Psi_r\left(\Lambda^r\right) = \int\limits_{\Psi^{-1}_r(A)}1_{\Psi^{-1}_r(\L_1[e])}d\Lambda^r=\sum\limits_{e'\in E^r,r(e')=e}\int\limits_{\Psi^{-1}_r(A)}1_{_1[e']}d\Lambda^r\\
&=&\sum\limits_{e'\in E^r,r(e')=e}\int\limits_{\Psi^{-1}_r(A)}\bar p^r_{e'}\circ F_r1_{\bar K^r_{i(e')}\setminus K^r_{i(e')}}\circ F_r1_{T^r_{i(e')}}d\Lambda^r\\
&=&\sum\limits_{e'\in E^r,r(e')=e}\sum\limits_{e''\in E^r,i(e'')=i(e')}\int\limits_{\Psi^{-1}_r(A)}\bar p_{e}\circ F_r1_{T^r_{i(e')}}\partial p^r_{e''}\circ F_rd\Lambda^r\\
&=&\sum\limits_{e'\in E^r,r(e')=e}\sum\limits_{e''\in E^r,i(e'')=i(e')}\int\limits_{\Psi^{-1}_r(A)}\bar p_{e}\circ F_r1_{_1[e'']} d\Lambda^r.
\end{eqnarray*}
Now, observe that, since $\M$ is positive and $\Lambda^r(\Sigma^r_G)=1$, by Lemma \ref{prm} (ii),
\[\sum\limits_{e'\in E^r,r(e')=e}\sum\limits_{e''\in E^r,i(e'')=i(e')}1_{_1[e'']} =\sum\limits_{e^r\in E^r,i(r(e^r))=i(e)}1_{_1[e^r]}=1_{T_{i(e)}}\circ\Psi_r\ \ \ \Lambda^r\mbox{-a.e.}\]
Hence, by Lemma \ref{esrl} (iii),
\begin{eqnarray*}
   \int\limits_A1_{_1[e]}d\Psi_r\left(\Lambda^r\right) = \int\limits_{\Psi^{-1}_r(A)}\bar p_{e}\circ F\circ\Psi_r1_{T_{i(e)}}\circ\Psi_r d\Lambda^r= \int\limits_{A}\bar p_{e}\circ F1_{T_{i(e)}} d\Psi_r\left(\Lambda^r\right).
\end{eqnarray*}
Therefore, since the class of sets of the form $_{n}[e_n,...,e_0]$ generates $\mathcal{F}$, is $\cap$-stable and covers $\bar\Sigma$, 
\[ E_{\Psi_r\left(\Lambda^r\right)}\left(1_{_1[e]}|\mathcal{F}\right) =\bar p_{e}\circ F1_{T_{i(e)}} \ \ \ \Psi_r\left(\Lambda^r\right)\mbox{-a.e.}\]
This completes the proof.
\hfill$\Box$

\begin{Proposition}\label{crp}
   Suppose $\mathcal{M}$ is uniformly continuous and positive such that  $D$ and $F|_D$ do not depend on the choice of $x_i$'s. Let $\M^r$ be a countable refinement of $\M$ which has a dominating Markov chain. Then $\M^r$ is consistent if $\M$ is consistent.
\end{Proposition}
{\it Proof.}  
 Let $\Lambda^r\in\tilde E(\M^r)$ and $f\in\mathcal{L}^B(K)$ be bounded. By Lemma \ref{rcl}, $\Psi_r(\Lambda^r)\in\tilde E(\M)$. Hence, since $\M$ is consistent, by Lemma \ref{esrl} (iii),
\begin{eqnarray*}
   \int UfdF_r\left(\Lambda^r\right) = \int UfdF\left(\Psi_r\left(\Lambda^r\right)\right) =\int fdF\left(\Psi_r\left(\Lambda^r\right)\right)= \int fdF_r\left(\Lambda^r\right).
\end{eqnarray*}
Thus, $F_r(\Lambda^r)\in P(\M^r)$.
\hfill$\Box$

\subsubsection{A recurrence condition}

It is well known from the theory of discrete homogeneous Markov chains that  the existence of an invariant probability measure requires the {\it positive recurrence} of the process. The following condition serves the same purpose for our generalization (see \cite{Wer13} for more elaboration on that). The next lemma gives then a sufficient condition for it in terms of the dominating Markov chain. (It was pointed out by an anonymous reviewer of \cite{Wer13} that the condition might be related to the notion of positive recurrence for countable Markov shifts introduced by O. Sarig, see \cite{Sa1} and \cite{Sa2}.)

\begin{Definition}
For $x\in K$, let $(\alpha^x_n)_{n\in\mathbb{N}}$ denote the sequence of probability measures on $N$ given by 
      \[\alpha^x_n(\{j\}):=\frac{1}{n}\sum\limits_{k=1}^{n}{U^*}^k\delta_{x}\left(K_j\right)\ \ \ \mbox{ for all }j\in N\mbox{ and }n\in\mathbb{N}.\]
\end{Definition}

\begin{Condition}\label{dtc}
There exists $x_0\in K$ such that $(\alpha^{x_0}_n)_{n\in\mathbb{N}}$
is uniformly tight, i.e. for every $\epsilon>0$ the exists a finite $V\subset N$ such that $\alpha^{x_0}_n(N\setminus V)<\epsilon$ for all $n\in\N$.
\end{Condition}

\begin{Lemma}\label{prl} 
    Condition \ref{dtc} is satisfied for all $x_0\in K$ if $\mathcal{M}$ has a  dominating Markov chain such that $c<\infty$.
\end{Lemma}
{\it Proof.}  
By the hypothesis,  $\xi_i<\infty$  for all $i\in N$, as in  \eqref{dmcc}.  Let  $k>0$.   Then for every $x\in K$  and $j\in N$,
\begin{eqnarray*}
{U^*}^k\delta_x\left(K_{j}\right)
&=&\sum\limits_{e_{-k},...,e_{-2}\in E}P^{-k}_{x}\left(\O_{-k}[e_{-k},..,e_{-2}]\right)\\
&&\times\sum\limits_{e\in E, t(e)=j}p_e\left( w_{\sigma_{-2}}\circ ... \circ w_{\sigma_{-k}}x\right)\\
&\leq&\sum\limits_{e_{-k},...,e_{-2}\in E}P^{-k}_{x}\left(\O_{-k}[e_{-k},..,e_{-2}]\right)\\
&&\times\sum\limits_{e\in E, i(e)=t(e_{-2}),t(e)=j}\sup\limits_{x\in K_{t(e_{-2})}}p_e(x)\\
&=&\sum\limits_{e_{-k},...,e_{-2}\in E}P^{-k}_{x}\left(\O_{-k}[e_{-k},..,e_{-2}]\right)\xi_{t(e_{-2})}q_{t(e_{-2})j}\\
&=&\sum\limits_{i\in N}\sum\limits_{e_{-k},...,e_{-2}\in E, t(e_{-2})=i}P^{-k}_{x}\left(\O_{-k}[e_{-k},..,e_{-2}]\right)\xi_{i}q_{ij}\\
&=&\sum\limits_{i\in N}{U^*}^{k-1}\delta_x\left(K_{i}\right)\xi_{i}q_{ij}
 \end{eqnarray*}  
where $(q_{ij})_{i,j\in N}$ is the transition matrix of the dominating Markov chain. Hence,
\begin{eqnarray*}
\alpha^x_{n+1}(\{j\})&\leq&\frac{1}{n+1}\sum\limits_{i\in N}\sum\limits_{k=1}^{n+1}{U^*}^{k-1}\delta_x\left(K_{i}\right)\xi_{i}q_{ij}\\
&\leq&\frac{1}{n+1}\sum\limits_{i\in N}\delta_x\left(K_{i}\right)\xi_{i}q_{ij}+\frac{n}{n+1}\sum\limits_{i\in N}\alpha^x_{n}(\{i\})\xi_{i}q_{ij}\\
&\leq&\frac{1}{n+1}\sum\limits_{i\in N}\delta_x\left(K_{i}\right)\xi_{i}q_{ij}+\sup\limits_{i\in N}\xi_{i}q_{ij}
 \end{eqnarray*} 
for all $n\in \N\cup\{0\}$ and $j\in N$. 

 Fix $x_0\in K$.  Let $i_0\in N$ such that $x_0\in K_{i_0}$. Let $\epsilon>0$. By the hypothesis, there exists a finite $V_\epsilon\subset N$ such that
 $\sum_{j\in N\setminus V_\epsilon}\sup_{i\in N}\xi_iq_{ij}<\epsilon/2$. Let $n_0\in\mathbb{N}$ such that $\xi_{i_0}q_{i_0j}/n_0<\epsilon/2$.    Then 
\[\alpha^{x_0}_n(N\setminus V_\epsilon)\leq\frac{\xi_{i_0}q_{i_0j}}{n}+\sum\limits_{j\in N\setminus V_\epsilon}\sup\limits_{i\in N}\xi_iq_{ij}<\frac{\epsilon}{2}+\frac{\epsilon}{2}=\epsilon\] for all $n\geq n_0$. 
Let $V_{n_0}\subset N$ be finite such that $\alpha^{x_0}_n(N\setminus V_{n_0})<\epsilon$ for all $1\leq n<n_0$. Then
\[\alpha^{x_0}_n\left(N\setminus (V_\epsilon\cup V_{n_0})\right)<\epsilon\]
for all $n\in \N$.
\hfill$\Box$

\subsection{Contractive, uniformly continuous Markov system}
In this subsection, we are going to apply the theory developed so far to the case when $\mathcal{M}$ is contractive.

Now, set \[L(x) := \sum\limits_{j\in N}d\left(x,x_j\right)1_{K_j}(x)\ \ \ \mbox{ for all }x\in K,\]
    \begin{equation}\label{cgc}
       b:= \sup\limits_{i\in N}\sup\limits_{x\in K_i}\sum\limits_{e\in E,\ i(e) = i}p_e(x)d\left( w_e(x_{i(e)}), x_{t(e)}\right),
    \end{equation}  
\[C(x):=L(x)+\frac{b}{1-a}\ \ \ \mbox{ for all }x\in K,\]
where $0<a<1$ is a contraction rate of the Markov system, and let us abbreviate $X_m(\sigma):= w_{\sigma_0}\circ...\circ w_{\sigma_{m}}(x_{i(\sigma_{m})})$ for all $\sigma\in\bar\Sigma$ and $m\leq 0$.

\begin{Lemma}\label{accl}
    Suppose $\mathcal{M}$ is contractive  with a contraction rate $0<a<1$.  \\
   (i) $UL \leq aL + b$.\\
   (ii)
\begin{equation}\label{ob}
    U^nL \leq a^nL+ \frac{b}{1-a}\ \ \ \mbox{ for all  } n\geq 0.
\end{equation}
(iii) For every $x\in K$, $m\leq 0$ and $n\geq 0$,
   \begin{equation}\label{si}
       \int d\left(X_m(\sigma),w_{\sigma_0}\circ ...\circ w_{\sigma_{m-n}}(x)\right) dP_x^{m-n}(\sigma)\leq a^{-m+1}C(x)
  \end{equation}
  and
\begin{eqnarray}\label{asC}
      \int d\left(X_m, X_{m-1}\right) dP_x^{m-n}\leq a^{-m+1} 2C(x).
\end{eqnarray}
\end{Lemma}
{\it Proof.}  
(i) Let $x\in K_i$ for some $i\in N$. Then
\begin{eqnarray*}
 UL(x)&=& \sum\limits_{e\in E}p_e(x)\sum\limits_{j\in N}d\left(w_e(x), x_j\right)1_{K_j}(w_ex)\\
      &\leq& \sum\limits_{e\in E,\ i(e) = i}p_e(x)d(w_e(x), w_e(x_i)) +  \sum\limits_{e\in E,\ i(e) = i}p_e(x)d\left(w_e(x_i), x_{t(e)}\right)\\
      &\leq& aL(x) + b.
\end{eqnarray*}

(ii) \eqref{ob} follows immediately from (i). 

(iii)  Clearly, the inequality is true if $b=\infty$. Now, suppose $b<\infty$. By  the contraction condition and \eqref{ob},
\begin{eqnarray*}
&&\int d\left(X_m(\sigma),w_{\sigma_0}\circ ...\circ w_{\sigma_{m}}\circ...\circ w_{\sigma_{m-n}}(x)\right) dP_x^{m-n}(\sigma)\\
&\leq&a^{-m+1}\sum\limits_{e_{m-n},...,e_{m-1}\in E}p_{e_{m-n}}(x)...p_{e_{m-1}}(w_{e_{m-2}}\circ ...\circ w_{e_{m-n}}x)\\
&&\ \ \ \ \ \ \ \ \ \ \ \ \ \ \ \ \ \ \ \ \ \ \times d\left(x_{t(e_{m-1})}, w_{e_{m-1}}\circ...\circ w_{e_{m-n}}(x)\right) \\
&=&a^{-m+1}U^{n-1}\left(\sum\limits_{e_{m-1}\in E}p_{e_{m-1}}d\left(x_{t(e_{m-1})}, w_{e_{m-1}}\right)\right)(x)\\
&\leq&a^{-m+1}\left(b+aU^{n-1}L(x)\right)\\
&\leq&a^{-m+1} C(x).
\end{eqnarray*}
This proves \eqref{si},  and \eqref{asC} follows by the triangle inequality.
\hfill$\Box$

\begin{Lemma}\label{pcl}
    Suppose $\mathcal{M}$ is contractive  with a contraction rate $0<a<1$. Then
\[\int L d\mu\leq\frac{b}{1-a}\ \ \ \mbox{ for all }\mu\in P(\mathcal{M}).\]
\end{Lemma}
{\it Proof.}  
Let $\mu\in P(\mathcal{M})$. Clearly, the inequality is true if $b=\infty$. Now, suppose $b<\infty$. For  $f\in\mathcal{L}^B(K)$ and $k\in\mathbb{N}$, set $f\wedge k:=\min\{f, k\}$. Observe that $U(f\wedge k)\leq k\wedge U(f)$. Hence, by induction, $U^n(f\wedge k)\leq k\wedge U^n(f)$ for all $n\in\mathbb{N}$. Therefore, by \eqref{ob}, 
\[\int L\wedge k d\mu = \int U^n(L\wedge k)d\mu\leq\int k\wedge\left(a^nL + \frac{b}{1-a}\right)d\mu \]
for all $n\in\mathbb{N}$. Since, for every $k\geq 0$, the functions $k\wedge(a^nL + b/(1-a))$ are integrable and converge monotonously to $k\wedge( b/(1-a))$, as $n\to\infty$, by the Monotone Convergence Theorem, 
\[\int L\wedge k d\mu \leq k\wedge\frac{b}{1-a}\]
for all $k\geq 0$. Applying the Monotone Convergence Theorem again, as $k\to\infty$, implies the assertion.
\hfill$\Box$

\subsubsection{Main theorem}

\begin{Theorem}\label{et} Suppose $\mathcal{M}$  is contractive with a contraction rate $0<a<1$ and uniformly continuous with $b<\infty$.
   Then the following holds true.\\
  (i)  If $\mathcal{M}$  has a dominating Markov chain, and Condition \ref{dtc} is satisfied for $x_0\in K$, then there exists $M\in\tilde E(\mathcal{M})$ such that for every $i\in N$ with $M(T_i)>0$ and $n\in\mathbb{N}$ there exist $k\geq n$ and a path $(e_1,...,e_k)$ with $t(e_k)=i$ such that $P_{x_0}(\L_1[e_1,...,e_k]^+)>0$. \\
(ii) $E(\mathcal{M})\subset\Phi(P(\mathcal{M}))\subset\tilde E(\mathcal{M})$ and $F(\Phi(\mu)) = \mu$ for all $\mu\in  P(\mathcal{M})$. \\
(iii) There exists a sequence of Borel sets $Q_1\subset
 Q_2\subset...\subset\Sigma_G$ with $\sum_{k\geq n} \Phi(\mu)(\bar\Sigma\setminus Q_k)\leq 1/(1-\sqrt{a})a^{n/2}$ for all  $\mu\in P(\mathcal{M})$ and $n\in\mathbb{N}$ such that  for each $k\in\mathbb{N}$ 
\[  d(F(\sigma),F(\sigma'))\leq \frac{8b}{(1-\sqrt{a})(1-a)}d'(\sigma,\sigma')^{\frac{\log\sqrt{a}}{\log(1/2)}}\]
whenever $\sigma,\sigma'\in Q_k\mbox{ with }d'(\sigma,\sigma')\leq(1/2)^{k+1}$, i.e. $F|_{Q_k}$ is locally H\"{o}lder-continuous with the same H\"{o}lder-constants for all $k\in\mathbb{N}$.\\
\end{Theorem}

{\it Proof.}  
 Set 
\[\phi^n_{m-1}:=\frac{1}{n}\sum_{k=1}^nP_{x_0}^{m-k}\]
for all $m\in\mathbb{Z}$ and $n\geq 1$. Then each $\phi^n_{m}$ is clearly a measure on $\mathcal{A}_m$.   Recall that $\bar\Sigma^+$ is a compact metrizable space.  Clearly, the set of all pre-images of cylinder sets in $\mathcal{B}(\bar\Sigma^+)$ under $\pi$ is exactly the set of all cylinder sets in $\mathcal{A}_{1}$.  Therefore, since both $\sigma$-algebras are generated by their cylinder sets, $\pi^{-1}(\mathcal{B}(\bar\Sigma^+)) =\mathcal{A}_{1}$, i.e. the induced set map $\tilde\pi^{-1}:\mathcal{B}(\bar\Sigma^+)\longrightarrow\mathcal{A}_{1}$ is bijective. Since the set of all Borel probability measures on $\bar\Sigma^+$ is  sequentially compact in the weakly-star topology, there exists a subsequence $(\pi(\phi^{n_k}_{1}))_{k\in\mathbb{N}}$ and a probability measure $\phi^+$ on $\mathcal{B}(\bar\Sigma^+)$ such that $\pi(\phi^{n_k}_{1})$ converges to $\phi^+$ weak-star as $k\to\infty$. Observe that, by the definition of  $\phi^{n}_{1}$'s, $\phi^+$ is invariant with respect to the left shift map on $\bar\Sigma^+$, as the shift maps commute with $\pi$. Set 
\[M(\pi^{-1}(B)):= \phi^+(B)\mbox { for all }B\in\mathcal{B}(\bar\Sigma^+).\]
Then this defines a shift-invariant measure $M$ on $\mathcal{A}_{1}$. Furthermore, by the shift-invariance of $M$, this gives consistent measures on $\mathcal{A}_m$ for all $m\in\mathbb{Z}$, by $M(S^{m-1}A)$,  which we  will also denote by $M$. In particular, $M$ defines consistent measures on all finite dimensional sub-$\sigma$-algebras of the product $\sigma$-algebra on $\bar\Sigma$. 
 Let $e_1,...,e_n\in E$. Observe that the image  of $\O  _{1}[e_1,...,e_n]\subset\bar\Sigma$ under $\pi$  is a cylinder set which is open and closed in $\bar\Sigma^+$.  Observe that $\phi^{n}_m = \phi^{n}_{1}\circ S^{m-1}$ for all $m\leq 1$ and $n\geq 1$. Therefore,
 \begin{eqnarray}\label{wcc}
     &&\lim\limits_{k\to\infty}\phi^{n_k}_m( _{m}[e_1,...,e_n] ) = \lim\limits_{k\to\infty}\phi^{n_k}_{1}\left(\O _{1}[e_1,...,e_n] \right) \nonumber\\
 &=&      \lim\limits_{k\to\infty}\phi^{n_k}_{1}\circ{\pi}^{-1}\left(\pi (\O _{1}[e_1,...,e_n])  \right)
= \phi^+\left(\pi (\O _{1}[e_1,...,e_n])  \right)\nonumber\\
& =&  M( _{m}[e_1,...,e_n] )
\end{eqnarray}
for all  $m\leq 1$. Furthermore, observe that
 \begin{eqnarray}\label{oswc}
\liminf\limits_{k\to\infty}\phi^{n_k}_m( O ) \geq M( O )
\end{eqnarray}
for every open set $O\in\mathcal{A}_{m}$ and $m\leq -1$, as $\pi$ is open.

 Now, for  every $m\leq 0$, $k\geq 1$ and  finite $C\subset E$,
\begin{eqnarray*}
&&P_{x_0}^{m-k}\left( \bigcup\limits_{e\in\bar E\setminus C} {_m[e]}\right)\\
&=&\sum\limits_{j\in N}\sum\limits_{e\in E\setminus C, i(e)=j}\int p_e\circ w_{\sigma_{m-1}}\circ ... \circ w_{\sigma_{m-k}}(x_0)P_{x_0}^{m-k}   (\sigma)\\
&=&\sum\limits_{j\in N}\sum\limits_{e\in E\setminus C, i(e)=j}U^k(p_e)(x_0).
 \end{eqnarray*}  
Let $\epsilon>0$. By the hypothesis, there exists a finite $V_\epsilon\subset N$ such that $\alpha^{x_0}_n(N\setminus V_\epsilon)<\epsilon/2$ for all $n\in\mathbb{N}$. Thus, by \eqref{oswc}, as $ \bigcup_{e\in\bar E\setminus C} {_m[e]}$ is open, 
\begin{eqnarray*}
    M\left( \bigcup\limits_{e\in\bar E\setminus C} {_m[e]}\right)&\leq&\liminf\limits_{k\to\infty}\phi^{n_k}_{m}\left( \bigcup\limits_{e\in\bar E\setminus C} {_m[e]}\right)\\
&\leq&\limsup\limits_{k\to\infty}\sum\limits_{j\in N}\sum\limits_{e\in E\setminus C, i(e)=j}\frac{1}{n_k}\sum\limits_{t=1}^{n_k}U^t(p_e)(x_0)\\
&\leq&\sum\limits_{j\in V_\epsilon}\sum\limits_{e\in E\setminus C, i(e)=j}\sup\limits_{x\in K_j}p_e(x) + \limsup\limits_{k\to\infty}\sum\limits_{j\in N\setminus V_\epsilon}\alpha^{x_0}_{n_k}(\{j\})\\
&\leq&\sum\limits_{j\in V_\epsilon}\sum\limits_{e\in E\setminus C, i(e)=j}\sup\limits_{x\in K_j}p_e(x) + \frac{\epsilon}{2}.
 \end{eqnarray*} 
Therefore,  by the hypothesis,  there exists a finite $C\subset E$ such that
\begin{equation}\label{caci}
    M\left( \bigcup\limits_{e\in\bar E\setminus C} {_m[e]}\right)\leq \epsilon.
 \end{equation} 
This means that each one-dimensional measure $M$ has an {\it approximating compact class}. Therefore, $M$ extends  uniquely to a shift-invariant Borel probability measure on $\bar\Sigma$,  which we will also denote by  $M$, e.g. by Kolmogorov Consistency Theorem  \cite{Bog}. (As $\bar E$ is a Polish space, the existence of a compact approximating class is actually automatic. However,\eqref{caci} is still needed  for the next step.)

Now, set \[\Omega_\infty:=\left\{\sigma\in\bar\Sigma|\mbox{ there exists }m\in\mathbb{Z}\mbox{ s.t. }\sigma_m=\infty\right\}.\]
By \eqref{caci}, for every $m\in\mathbb{Z}$ there exists a finite $C_m\subset E$ such that
\[M\left( \bigcup\limits_{e\in\bar E\setminus C_m} {_m[e]}\right)\leq\frac{\epsilon}{4}2^{-|m|}.\]
As $\Omega_\infty\subset \bigcup_{m\in\mathbb{Z}}\bigcup_{e\in\bar E\setminus C_m} {_m[e]}$, it follows that $M(\Omega_\infty)\leq 3/4\epsilon<\epsilon$. Since $\epsilon$ was arbitrary, we conclude that
\begin{equation}\label{zmoi}
 M\left(\Omega_\infty\right)=0.
 \end{equation}  

Next, we show that $\Sigma_G$ has the full measure. First, observe that every $\sigma\in\bar\Sigma\setminus\Sigma_G$ is either in $\Omega_\infty$, or there exists $e_1,...,e_n\in E$ such that $(e_1,...,e_n)$ is not a path and $\sigma\in\O _m[e_1,...,e_n]$ for some $m\in\mathbb{Z}$. By Remark \ref{ppr} and \eqref{wcc}, $M( _m[e_1,...,e_n])=0$. Hence, for every $\sigma\in\bar\Sigma\setminus(\Omega_\infty\cup\Sigma_G)$ there exists an open set $O_\sigma$ such that $\sigma\in O_\sigma$ and $M(O_\sigma) = 0$.  Choose open sets $O_\infty\subset\bar\Sigma$  and $O_G\subset\bar\Sigma$ such that $\Omega_\infty\subset O_\infty$, $\Sigma_G\subset O_G$, $M( O_\infty)<\epsilon/2$ and $M(O_G\setminus\Sigma_G)<\epsilon/2$. Then, by the compactness of $\bar\Sigma$, there exist finitely many $\sigma^1,...,\sigma^k\in\bar\Sigma\setminus(\Omega_\infty\cup\Sigma_G) $ such that
$\bar\Sigma=O_G\cup O_\infty\cup\bigcup_{i=1}^kO_{\sigma^i}$. Hence, $M(O_G\cup\ O_\infty)=1$, and therefore, $M(\Sigma_G)>1-\epsilon$. Since $\epsilon$ was arbitrary, we conclude that
\begin{equation}\label{fmps}
   M\left(\Sigma_G\right)=1.
\end{equation}

 Now, we are going to show that $M(D) = 1$.  For $x\in K$ and $m\leq 0$, set 
\[A_x^{m-1}:=\left\{\sigma\in\bar\Sigma\left|\  d\left(X_m(\sigma),X_{m-1}(\sigma)\right)>a^{\frac{-m+1}{2}}2C(x)\right.\right\}. \]
 Then,   by \eqref{asC}, 
\[P^{m-n}_x\left(A^{m-1}_x \right)\leq a^{\frac{-m+1}{2}}\]
for all $x\in K$, $m\leq 0$ and $n\geq 1$. Hence
\begin{equation}\label{phib}
 \phi^{n}_m\left(A^{m}_{x_0} \right)\leq a^{\frac{-m}{2}}
\end{equation}
for all $m\leq -1$ and $n\geq 1$. Therefore, by \eqref{zmoi} and \eqref{oswc}, as  $A^{m}_{x_0}\\ \cap\bigcup_{e_m,...,e_0\in E}\O _m[e_m,...,e_0]$ is  a countable  union of some open cylinder sets,
\[M\left(A^{m}_{x_0} \right) = M\left(A^{m}_{x_0}\setminus\Omega_\infty \right) \leq\liminf\limits_{k\to\infty}\phi^{n_k}_m\left( A^{m}_{x_0}\cap\bigcup_{e_m,...,e_0\in E}\O _m[e_m,...,e_0]\right)\leq a^{\frac{-m}{2}}\]
for all $m\leq -1$.

Now, set
\[A_{x_0}:=\bigcap\limits_{l\leq -1}\bigcup\limits_{m\leq l}A^{m}_{x_0}.\]
Then
\[M\left(A_{x_0}\right)\leq\sum\limits_{m\leq l}M\left(A^{m}_{x_0} \right)\leq \sum\limits_{m\leq l}a^{\frac{-m}{2}}\mbox{ for all }l\leq -1.\]
Hence
\[M\left(A_{x_0}\right) = 0.\]
Now, observer that for every $\sigma\in\Sigma_G\setminus A_{x_0}$, sequence $(X_m(\sigma))_{m\leq 0}$ is Cauchy. Hence, by the completeness of $(K, d)$, $\Sigma_G\setminus A_{x_0}\subset D$. Therefore,
\begin{equation}\label{fmd}
   M\left(D\right) = 1.
\end{equation}

Now, we are going to compute $E_M(1_{_1[e]}|\mathcal{F})$ for all $e\in E$.  Fix $e\in E$. First, observe that, for $x\in K$, $m\leq 0$, $n\geq 0$ and  $_m[e_m,...,e_0]\subset\bar\Sigma$,
 \begin{eqnarray}\label{dei}
     &&\int\limits_{_m[e_m,...,e_0]}1_{_1[e]}dP^{m-n}_x\nonumber\\
&=&\sum\limits_{e_{m-n},...,e_{m-1}\in E}P^{m-n}_x\left(_{m-n}[e_{m-n},...,e_{m-1},e_m,...,e_0,e ]\right)\nonumber\\
&=&U^{n}\left(P^m_.( _m[e_m,...,e_0])p_e\circ w_{e_0}\circ ... \circ w_{e_m}\right)(x)\nonumber\\
&=&\int\int\limits_{_m[e_m,...,e_0]}p_e\circ w_{\sigma_0}\circ ... \circ w_{\sigma_m}(y)dP^{m}_y(\sigma)d{U^*}^n\delta_x(y)\nonumber\\
&=&\int\limits_{_m[e_m,...,e_0]}p_e\circ X_m(\sigma)dP^{m-n}_x(\sigma) + r_{mn,x}\left( _m[e_m,...,e_0]\right) 
\end{eqnarray}
where $r_{mn,x}$ is a signed measure on $\mathcal{F}_m$ given by
 \begin{eqnarray*}
    r_{mn,x}\left( A\right):=\int\int\limits_{A}\left(p_e\circ w_{\sigma_0}\circ ... \circ w_{\sigma_m}(y)-p_e\circ X_m(\sigma)\right)dP^{m}_y(\sigma)d{U^*}^n\delta_x(y)
\end{eqnarray*}
for all $A\in \mathcal{F}_m$. Hence, as every member of $ \mathcal{F}_m$ can be written as a countable disjoint union of cylinder sets,
\begin{equation*}
      \int\limits_{A}1_{_1[e]}d\phi^n_{m} = \int\limits_{A}p_e\circ X_m d\phi^n_{m} + \frac{1}{n}\sum\limits_{k=1}^nr_{mk,x_0} (A)
\end{equation*}
for all $A\in \mathcal{F}_m$ and $m< 0$ and $n\in\mathbb{N}$. This implies, by \eqref{wcc} and \eqref{zmoi}, that
\begin{equation}\label{cei}
      \int\limits_{A}1_{_1[e]}dM = \int\limits_{A}p_e\circ X_m dM + \lim\limits_{k\to\infty}\frac{1}{n_k}\sum\limits_{k=1}^{n_k}r_{mk,x_0} (A)
\end{equation}
for all $A$ which are finite unions of cylinder sets from $\mathcal{F}_m$ and $m< 0$. 
Now, for $y\in K$ and $m\leq 0$, set
\[B_{m,y}:=\left\{\sigma\in\bar\Sigma|\ d(w_{\sigma_0}\circ ... \circ w_{\sigma_m}(y),X_m(\sigma))>a^{\frac{-m+1}{2}}C(y)\right\}\]
and 
\[\beta_e(t):=\sup\limits_{x,y\in K_{i(e)},d(x,y)\leq t}\left\{p_e(x) - p_e(y)\right\}\mbox{ for all }t\geq 0.\]
Then, by \eqref{si},
\[P^{m}_y\left(B_{m,y}\right)\leq a^{\frac{-m+1}{2}}\mbox{ for all }y\in K\mbox{ and } m\leq 0.\] 
Therefore, 
 \begin{eqnarray}\label{refa}
     \left|r_{mk,x_0} (A)\right| &\leq& a^{\frac{-m+1}{2}} +\int \beta_e\left(a^{\frac{-m+1}{2}}C(y)\right)d{U^*}^k\delta_{x_0}(y)
\end{eqnarray}
for all $A\in \mathcal{F}_m$, $m\leq 0$ and $k\in\mathbb{N}$.  
Set $\rho :=  b/(1-a)+L(x_0)$ and 
\[B (\alpha):=\bigcap\limits_{j\in N}\left(K\setminus K_j\right)\cup B_{\alpha}(x_j)\]
for all $\alpha\geq 0$. Then, by \eqref{ob},
\begin{eqnarray*}
    \rho&\geq& U^kL(x_0) \\
   &=& \int\sum\limits_{j\in N} d(w_{\sigma_k}\circ ... \circ w_{\sigma_1}(x_0), x_j) 1_{K_j}\circ w_{\sigma_k}\circ ... \circ w_{\sigma_1}(x_0)dP^1_{x_0}(\sigma)\\
             &\geq& \frac{2\rho}{\epsilon}\sum\limits_{j\in N}
            P^1_{x_0}\left(d(w_{\sigma_k}\circ ... \circ w_{\sigma_1}(x_0), x_j) >  \frac{2\rho}{\epsilon}\mbox{ and }w_{\sigma_k}\circ ... \circ w_{\sigma_1}(x_0)\in K_j\right)\\
&\geq&  \frac{2\rho}{\epsilon}P^1_{x_0}\left(w_{\sigma_k}\circ ... \circ w_{\sigma_1}(x_0)\in K\setminus B\left( \frac{2\rho}{\epsilon}\right)\right)
\end{eqnarray*}
 for all $k\geq 1$. Hence
 \begin{equation*}
       P^1_{x_0}\left(w_{\sigma_k}\circ ... \circ w_{\sigma_1}(x_0)\in K\setminus B\left( \frac{2\rho}{\epsilon}\right)\right)\leq\frac{\epsilon}{2}\mbox{ for all }k\geq 1.
\end{equation*}
 Thus
\begin{eqnarray*}
     {U^*}^k\delta_{x_0}\left(K\setminus B\left( \frac{2\rho}{\epsilon}\right)\right) &=& \int 1_{ K\setminus B\left( \frac{2\rho}{\epsilon}\right)}\circ w_{\sigma_k}\circ ... \circ w_{\sigma_1}(x_0) dP^1_{x_0}(\sigma)\\
&=& P^1_{x_0}\left(w_{\sigma_k}\circ ... \circ w_{\sigma_1}(x_0)\in K\setminus  B\left( \frac{2\rho}{\epsilon}\right)\right)\\
&\leq&\frac{\epsilon}{2}
\end{eqnarray*}
for all $k\in\mathbb{N}$. Now, set $C_\epsilon:=2\rho/\epsilon +b/(1-a)$. Then $C(y)\leq C_\epsilon$ for all $y\in B( 2\rho/\epsilon)$. Therefore, by \eqref{refa},
    \[ \left|r_{mk,x_0} (A)\right| \leq a^{\frac{-m+1}{2}} + \beta_e\left(a^{\frac{-m+1}{2}}C_\epsilon\right)+\frac{\epsilon}{2}\]
for all  $A\in \mathcal{F}_m$, $m\leq 0$ and $k\in\mathbb{N}$.  Thus, by \eqref{cei},
\begin{equation}\label{cma}
      \left|\int\limits_{A}1_{_1[e]}dM - \int\limits_{A}p_e\circ X_mdM\right|\leq a^{\frac{-m+1}{2}} + \beta_e\left(a^{\frac{-m+1}{2}}C_\epsilon\right)+\frac{\epsilon}{2}
\end{equation}
for all $A$ which are finite unions of cylinder sets from $\mathcal{F}_m$ and $m< 0$, and, since every member of $\mathcal{F}_m$ can be written as a countable union of cylinder sets, by  Lebesgue's Dominated  Convergence Theorem, it holds true for all $A\in\mathcal{F}_m$ and $m< 0$. That is
\begin{equation*}
      \left|\int\limits_{A}\left(E_M\left(1_{_1[e]}|\mathcal{F}_m\right) - p_e\circ X_m\right)dM\right|\leq a^{\frac{-m+1}{2}} + \beta_e\left(a^{\frac{-m+1}{2}}C_\epsilon\right)+\frac{\epsilon}{2}
\end{equation*}
for all $A\in \mathcal{F}_m$ and $m< 0$. Set $A_m^-:=\{\sigma\in\bar\Sigma|\ E\left(1_{_1[e]}|\mathcal{F}_m\right)(\sigma) \leq p_e\circ X_m(\sigma)\}$. Then
\begin{equation*}
      \int\limits_{A^-_m}\left|E_M\left(1_{_1[e]}|\mathcal{F}_m\right) - p_e\circ X_m\right|dM\leq a^{\frac{-m+1}{2}} + \beta_e\left(a^{\frac{-m+1}{2}}C_\epsilon\right)+\frac{\epsilon}{2} 
\end{equation*}
and, obviously, the same inequality holds true also with $\bar\Sigma\setminus A^-_m$ in place of $A^-_m$ for all $m< 0$. Hence
\begin{equation}\label{ceia}
      \int\limits\left|E_M\left(1_{_1[e]}|\mathcal{F}_m\right) - p_e\circ X_m\right|dM\leq 2a^{\frac{-m+1}{2}} +2 \beta_e\left(a^{\frac{-m+1}{2}}C_\epsilon\right)+\epsilon 
\end{equation}
for all $m< 0$. Let $\sigma\in D$. Observe that $X_m(\sigma)\in K_{t(\sigma_0)}$ for all $m\leq 0$. Therefore,  $\lim_{m\to-\infty}p_e\circ X_m(\sigma) = \bar p_e\circ F(\sigma)$ if $i(e) =t(\sigma_0)$. Otherwise,
$\lim_{m\to-\infty}p_e\circ X_m(\sigma) = 0$. Hence, since $M(D)=1$,
\[\lim_{m\to-\infty}p_e\circ X_m(\sigma) =  \bar p_e\circ F(\sigma)1_{K_{i(e)}}\left(x_{t(\sigma_0)}\right)\mbox{ for }M\mbox{-a.a. }\sigma\in\bar\Sigma. \]
Therefore, by Lebesgue's Dominated Convergence Theorem, $p_e\circ X_m(\sigma)$  converges to  $ \bar p_e\circ F(\sigma)1_{K_{i(e)}}(x_{t(\sigma_0)})$ in $\mathcal{L}^1(M)$.
  Therefore,  by the triangle inequality and \eqref{ceia}, as $E_M\left(1_{_1[e]}|\mathcal{F}_m\right)$ also converges to $E_M\left(1_{_1[e]}|\mathcal{F}\right)$ in $\mathcal{L}^1(M)$,
\begin{eqnarray*}
   \int\limits\left|E_M\left(1_{_1[e]}|\mathcal{F}\right)(\sigma) - \bar p_e\circ F(\sigma)1_{K_{i(e)}}\left(x_{t(\sigma_0)}\right)\right|dM(\sigma)\leq \epsilon.
\end{eqnarray*}
Since $\epsilon$ was arbitrary, we conclude that
\begin{equation*}
     E_M\left(1_{_1[e]}|\mathcal{F}\right)(\sigma) = \bar p_e\circ F(\sigma)1_{K_{i(e)}}\left(x_{t(\sigma_0)}\right)\ \ \ \mbox{ for }M\mbox{-a.a. }\sigma\in\bar\Sigma
\end{equation*}
for all $e\in E$.  Thus $M\in\tilde E(\mathcal{M})$.  

Now, let $i\in N$ such that $M(T_i)>0$. Since $M(\Sigma_G)=1$, there exists $e\in E$ with $i(e)=i$ such that $M(T_i\cap\L_1[e])>0$.  Then for every $m\leq 0$ there exists a finite union $A_m:=\bigcup_{e_m,...,e_0}\L _m[e_m,...,e_0]$ where each  $(e_m,...,e_0, e)$ is a path and $M(A_m\cap\L_1[e])>M(T_i\cap\L_1[e])/2$. Let $0<\epsilon<M(T_i\cap\L_1[e])$. Choose  $m_0<0$ such that $a^{(-m_0+1)/2} + \beta_e\left(a^{(-m_0+1)/2}C_\epsilon\right)+\epsilon/2< M(T_i\cap\L_1[e])/2$. Then, by \eqref{cma},
 \begin{equation*}
    0<M(A_{m_0}\cap\L_1[e])-a^{\frac{-m_0+1}{2}} -\beta_e\left(a^{\frac{-m_0+1}{2}}C_\epsilon\right)-\frac{\epsilon}{2}\leq \int\limits_{A_{m_0}}p_e\circ X_{m_0}dM.
\end{equation*}
Hence, there exists a path $(e_{m_0},...,e_0, e)$ such that
 \begin{eqnarray*}
   0&<&p_e\circ X_{m_0}M(\L _{m_0}[e_{m_0},...,e_0])\\
&\leq&\lim\limits_{k\to\infty}\frac{1}{n_k}\sum_{j=1}^{n_k}P_{x_0}^{m_0-j}(\L _{m_0}[e_{m_0},...,e_0]).
\end{eqnarray*}
Thus, for any $n\in\mathbb{N}$ there exists $j\geq n$ and a path $(e_{m_0-j},...e_{m_0-1},e_{m_0},...,e_0)$ with $t(e_0)=i$ such that
$P_{x_0}^{m_0-j}(\L _{m_0-j}[e_{m_0-j},...e_{m_0-1},e_{m_0},...,e_0])>0$.
This completes the proof of (i).

(ii) Inclusion $E(\mathcal{M})\subset\Phi(P(\mathcal{M}))$ follows from Lemma \ref{pesl} (i) and Proposition \ref{cml}. 

Now,  we show that  $ \Phi(P(\mathcal{M}))\subset\tilde E(\mathcal{M})$.  Set $G_m:=\{\sigma\in\bar\Sigma|\ (\sigma_m,...,\sigma_{|m|})\mbox{ is a path}\}$ for all $m<0$. Then,  by Remark \ref{ppr},  $\Phi(\mu)(G_m) = 1$ for all $m<0$. As  $\Sigma_G=\bigcap_{m<0}G_m$ and $G_{m-1}\subset G_m$ for all $m<0$, it follows that  
\begin{equation}\label{sftfm}
   \Phi(\mu)(\Sigma_G) = 1.
\end{equation}
 The integration of \eqref{asC} with respect to $\mu$ implies that
 \[\int d\left(X_m,X_{m-1}\right) d\Phi(\mu)\leq a^{-m+1} 2\int C(x) \mu(x)\]
for all $m\leq 0$. By Lemma \ref{pcl},  $\int C(x) \mu(x)\leq 2b/(1-a)$, therefore we can define
\[A_{m-1}:=\{\sigma\in\bar\Sigma|\ d\left(X_m(\sigma),X_{m-1}(\sigma)\right)>a^{\frac{-m+1}{2}}2C_1\}\]
with $C_1:= 2b/(1-a)$ for all $m\leq 0$. Then
\begin{equation}\label{iacp}
 \Phi(\mu)\left(A_{m} \right)\leq a^{\frac{-m}{2}}
\end{equation}
for all $m\leq -1$, and the same way as for \eqref{fmd}, this implies that 
\[\Phi(\mu)(D) = 1.\] 
Now, let $e\in E$. The integration of \eqref{dei} with respect to $\mu$ or the straightforward usage of the definition of $\Phi(\mu)$ gives
 \begin{eqnarray}\label{cepr}
     \int\limits_{_m[e_m,...,e_0]}1_{_1[e]}d\Phi(\mu)
=\int\limits_{_m[e_m,...,e_0]}p_e\circ X_md\Phi(\mu) + r_{m}\left( _m[e_m,...,e_0]\right) 
\end{eqnarray}
for all $_m[e_m,...,e_0]\in \mathcal{F}_m$ where $r_{m}$ is  given by
 \begin{eqnarray*}
    r_{m}\left( A\right):=\int\int\limits_{A}\left(p_e\circ w_{\sigma_0}\circ ... \circ w_{\sigma_m}(y)-p_e\circ X_m(\sigma)\right)dP^{m}_y(\sigma)d\mu(y)
\end{eqnarray*}
for all $A\in \mathcal{F}_m$. Furthermore, the integration of \eqref{si} gives
\begin{equation}\label{ide}
 \int \int d\left(w_{\sigma_0}\circ ...\circ w_{\sigma_{m}}(y), X_m(\sigma)\right) dP^{m}_y(\sigma)d\mu(y)\leq a^{-m+1}C_1
\end{equation}
  for all $m\leq 0$. Hence, the same way as for \eqref{refa}, it follows that
 \begin{eqnarray}\label{ire}
     \left|r_{m} (A)\right| &\leq& a^{\frac{-m+1}{2}} + \beta_e\left(a^{\frac{-m+1}{2}}C_1\right)
\end{eqnarray}
for all $A\in \mathcal{F}_m$. Therefore, as in the proof of (i), \eqref{cepr}  implies  that
 \begin{equation*}
      E_{\Phi(\mu)}\left(1_{_1[e]}|\mathcal{F}\right) = \bar p_e\circ F1_{T_{i(e)}}\ \ \ \Phi(\mu)\mbox{-a.e.} 
\end{equation*}
for all $e\in E$.  Thus $\Phi(\mu)\in\tilde E(\mathcal{M})$, as desired.

Next, we show that $F(\Phi(\mu)) = \mu$. It is sufficient to show that the measures agree on all  bounded uniformly continuous non-negative functions on $K$ (as this set of functions is closed under multiplication and generates Borel $\sigma$-algebra). Let  $f\in\mathcal{L}^B(K)$ be  bounded and uniformly continuous. Observe that, for each $m\leq 0$,
\begin{eqnarray}\label{imc}
    \int fd\mu&=& \int U^{-m+1}(f)d\mu\nonumber\\
      &=&  \int \int f\circ w_{\sigma_0}\circ ...\circ w_{\sigma_{m}}(x)dP^m_x(\sigma)d\mu(x)\nonumber\\
      &=& \int f\circ X_md\Phi(\mu) + R_{m}
\end{eqnarray}
where  
 \begin{eqnarray*}
    R_{m}:=\int\int\left(f\circ w_{\sigma_0}\circ ... \circ w_{\sigma_m}(y)-f\circ X_m(\sigma)\right)dP^{m}_y(\sigma)d\mu(y).
\end{eqnarray*}
Since $f$ is uniformly continuous and bounded, one sees, by the contraction condition, the same way as for \eqref{ire},
that $| R_{m}|\to 0$ as $m\to-\infty$. Therefore, since $\Phi(\mu)(D) = 1$, by the already shown, \eqref{imc} implies by Lebesgue's Dominated Convergence Theorem  that
\[\int fd\mu  = \int f\circ Fd\Phi(\mu) =\int fdF(\Phi(\mu)), \]
as desired.

(iii) Now, set
\[Q_k:=\bigcap\limits_{m<-k}\bar\Sigma\setminus A_m\]
for all $k\in\mathbb{N}$. By  \eqref{iacp}, $\sum_{k\geq n}\Phi(\mu)(\bar\Sigma\setminus Q_k) \leq 1/(1-\sqrt{a})a^{n/2}$ for all $\mu\in P(\mathcal{M})$ and  $n\in\mathbb{N}$. Since $\Phi(\mu)(D)=1$ for all $\mu\in P(\mathcal{M})$,  we can assume $Q_k\subset D$ for all $k$.  The proof that $F|_{Q_k}$ is locally H\"{o}lder continuous is the same as that of Lemma 3 (iii) in \cite{Wer3}. We give it for completeness here. Let $\sigma,\sigma'\in
Q_l$ for some $l\in\mathbb{N}$. Then, by the triangle inequality,
\[d\left(X_m(\sigma),X_{m-k}(\sigma)\right)\leq\sum\limits_{i\leq m}
a^{\frac{-i+1}{2}}2C_1=2C_1\frac{1}{1-\sqrt{a}}\ a^{\frac{-m+1}{2}}\mbox{ for all }m\leq -l,\ k\geq 1.\] Hence
\[d\left(X_m(\sigma),F(\sigma)\right)\leq 2C_1\frac{1}{1-\sqrt{a}}\ a^{\frac{-m+1}{2}}\mbox{ for all }m\leq -l.\]
 The same way,
 \[d\left(X_m(\sigma'),F(\sigma')\right)\leq 2C_1\frac{1}{1-\sqrt{a}}\ a^{\frac{-m+1}{2}}\mbox{ for all }m\leq -l.\] 
Now, let  $d'(\sigma,\sigma')=(1/2)^{-m+1}$ for some $m\leq -l$. Then
$X_m(\sigma')=X_m(\sigma)$.
Therefore,
\[d\left(F(\sigma),F(\sigma')\right)\leq\frac{4C_1}{1-\sqrt{a}}\ a^{\frac{-m+1}{2}}=\frac{8b}{(1-\sqrt{a})(1-a)}\ d'(\sigma,\sigma')^{\log\sqrt{a}/\log(1/2)}.\]
 This completes the proof of the theorem. 
\hfill$\Box$

\begin{Remark}
 An anonymous  reviewer pointed out that it might be possible to replace the compactification with the tightness argument for obtaining measure $M$ in the proof of Theorem  \ref{et} and restrict the consideration to $\Sigma_G$. Also, it might be possible to replace other mertizable compactness arguments in the article with the completeness and the separability arguments. However,  the author thinks that it would be short-sighted to discard the powerful information of working on a compact metrizable space from the set-up of the theory, in particular, because it gives the access to the results of the well-developed ergodic theory on the standard topological space, e.g. the author refers to it in the follow-up paper \cite{Wer13}, which uses many results from this paper.  Note that, in contrast to the mainstream thermodynamic formalism on countable Markov shifts, we do not have the continuity of the potential on  $\Sigma_G$, and therefore, we have nothing to lose by working on the standard topological space of the ergodic theory, and sometimes taking advantage of it. Furthermore, with not having the openness of the Markov partition, and therefore, working with possibly overlapping closures of the atoms of the partition, it seems to be unreasonable to restrict {\it a priori} the consideration only to $\Sigma_G$.

 Also, the reviewer found that 'The proof of theorem 5 could be simplified (and significantly shortened) by using tightness to prove the existence of an invariant measure, Borel-Cantelli to prove $M (D) = 1$
and Radon-Nikodym derivatives instead of conditional expectations.', though no proof was presented.
\end{Remark}

\begin{Corollary}\label{ndce} 
     Suppose $\mathcal{M}$  is contractive, uniformly continuous and  non-degenerate such that $b<\infty$. Then the following holds true.\\
(i) Suppose $\mathcal{M}$  has a dominating Markov chain. Then $E(\mathcal{M})$ is not empty if and only if Condition \ref{dtc} is satisfied.\\
(ii)   $\Phi$ is the inverse of $F: E(\mathcal{M})\longrightarrow P(\mathcal{M})$. (Thus, the latter does not depend on the choice of $x_i\in K_i$ for all $i\in N$  as long $b$ remains finite.)
\end{Corollary}
{\it Proof.}  
(i) The 'only if' part follows from Lemma \ref{pesl} (iv). The 'if' part follows from Theorem \ref{et} (i) and Theorem \ref{cndc}.

(ii) The assertion  follows by Theorem \ref{et} (iii), Theorem \ref{cndc} and Lemma \ref{pesl} (i).
\hfill$\Box$

\subsubsection{Invariant measures}

\begin{Corollary}\label{eimc} 
     Suppose $\mathcal{M}$  is contractive and uniformly continuous such that $b<\infty$.
   Then the following holds true.\\
   (i) Suppose $\mathcal{M}$ is non-degenerate and has a dominating Markov chain. Then $P(\mathcal{M})$ is not empty if and only if Condition \ref{dtc} is satisfied.\\
   (ii) Suppose $\mathcal{M}$ is consistent, has a dominating Markov chain and satisfies Condition \ref{dtc}. Then $P(\mathcal{M})$ is not empty.\\
   (iii)   Every $\mu \in P(\mathcal{M})$ is  tight.\\
   (iv) Suppose  sets  $\{e\in E|\ i(e)=j\}$ and $\{e\in E|\ t(e)=j\}$ are finite for all $j\in N$. Then for every finite $c\subset N$ and $\epsilon>0$ there exists a compact $C\subset K$ such that
\[\mu\left(C\right)>\mu\left(\bigcup\limits_{i\in c}K_i\right)-\epsilon\ \ \ \mbox{ for all }\mu\in P(\M).\]
  (v) If $E$ is finite, then $P(\M)$ is uniformly tight.
\end{Corollary}
{\it Proof.}  
(i) The assertion follows by Corollary \ref{ndce}  and Proposition \ref{cml}.

(ii)  The assertion follows by Theorem \ref{et} (i).

(iii) Let $\epsilon>0$. By Theorem \ref{et} (iii), there exists a Borel set $Q\subset\bar\Sigma$  with $Q\subset\Sigma_G$ such that $\Phi(\mu)(Q)>1-\epsilon/2$ for all $\mu\in P(\mathcal{M})$ and $F|_{Q}$ is uniformly continuous with respect to $d'$. Let $A\subset Q$ be compact in $(\Sigma_G, d')$. Then $C:=F(A)$ is compact, and,  by Theorem \ref{et} (ii),
\begin{equation}\label{ti}
   \mu\left( C\right)=\Phi(\mu)\left(F^{-1}\left( C\right)\right)\geq\Phi(\mu)\left(A\right)
\end{equation}
for all $\mu\in P(\M)$. Using the facts that  $\bar\Sigma$ compact, and every open ball in $(\Sigma_G, d')$ is contained in a cylinder set $_{-n}[e_{-n}, ...,e_n]$ with $e_{-n}, ...,e_n\in E$ such that $(e_{-n}, ...,e_n)$ is a path, which is closed in $\bar\Sigma$,  one easily checks that  $(\Sigma_G, d')$ is complete, and therefore, it is Polish. Let $\mu\in P(\M)$. Then, since $d'$ generates exactly the topology on  $\Sigma_G$ which is induced from $\bar\Sigma$, the restriction of  $\Phi(\mu)$ on $\Sigma_G$ is a Radon measure, and, by Theorem 5 (ii),  $\Phi(\mu)(\Sigma_G)=1$. Hence, there exists $C'\subset Q$ which is compact in $(\Sigma_G, d')$ such that  $\Phi(\mu)(Q\setminus C')<\epsilon/2$.   Therefore, by \eqref{ti}, setting $A:=C'$ implies that
  \[\mu\left(C\right)\geq\Phi(\mu)\left(C'\right)=\Phi(\mu)\left(Q\right)-\Phi(\mu)(Q\setminus C')>1-\frac{\epsilon}{2}-\frac{\epsilon}{2}=1-\epsilon.\] 
That is, $\mu$ is tight.

(iv)  Let   $e\in E$. Since $\bar\Sigma\setminus\L_0[e]$ is open in $\bar\Sigma$, and $d'$ generates exactly the topology on  $\Sigma_G$ which is induced from $\bar\Sigma$, $\Sigma_G\setminus\L_0[e]=\Sigma_G\cap\bar\Sigma\setminus\L_0[e]$ is open in $(\Sigma_G, d')$. Therefore,
$\Sigma_G\cap\L_0[e]=\Sigma_G\setminus(\Sigma_G\setminus\L_0[e])$ is closed in $(\Sigma_G, d')$. Now, we show that $\Sigma_G\cap\L_0[e]$ is compact in $(\Sigma_G, d')$. Let $(\sigma^k)_{k\in\N}\subset\Sigma_G\cap\L_0[e]$. By the hypothesis, for every $n\in\N$, $\Sigma_G\cap\L_0[e]$ can be written as a finite union of sets $\Sigma_G\cap\L_{-n}[e_{-n}, ...,e_n]$ where each $(e_{-n}, ...,e_n)$ is a path. Since the sets are open balls in $(\Sigma_G, d')$, choosing $\sigma^{k_n}$ from the ball containing infinite number of $\sigma^k$'s for each $n\in\N$ gives a Cauchy subsequence, which, by the completeness of  $(\Sigma_G, d')$, converges to some $\sigma\in\Sigma_G$. Since $\Sigma_G\cap\L_0[e]$ is closed in $(\Sigma_G, d')$, $\sigma\in\Sigma_G\cap\L_0[e]$. Thus $\Sigma_G\cap\L_0[e]$ is compact in $(\Sigma_G, d')$.
Now, let $\tilde Q$ be the closure of $Q$ in $(\Sigma_G, d')$, and $\tilde F$ be the continuous extension of $F|_{Q}$ on $\tilde Q$ in $(\Sigma_G, d')$. Then  $\tilde Q\cap\L_0[e]=\tilde Q\cap(\Sigma_G\cap\L_0[e])$ is compact in  $(\Sigma_G, d')$. Let $E_c:=\{e\in E|\ i(e)\in c\}$. Then, by the hypothesis, $E_c$ is finite, and therefore, 
\[B:=\tilde Q\cap\bigcup\limits_{e\in E_c}\L_0[e]\]
is compact in $(\Sigma_G, d')$. Hence, $\tilde C:=\tilde F(B)$ is compact in $(K,d)$, and, by Theorem \ref{et} (ii), for every $\mu\in P(\M)$,
\begin{eqnarray*}
\mu\left(\tilde C\right)&=&\Phi(\mu)\left(F^{-1}\left(\tilde C\right)\right)=\Phi(\mu)\left(F^{-1}\left(\tilde F\left(\tilde Q\cap\bigcup\limits_{e\in E_c}\L_0[e]\right)\right)\right)\\
&\geq&\Phi(\mu)\left(Q\cap\bigcup\limits_{e\in E_c}\L_0[e]\right)>1-\frac{\epsilon}{2} - \Phi(\mu)\left(Q\setminus\bigcup\limits_{e\in E_c}\L_0[e]\right)\\
&\geq& 1-\frac{\epsilon}{2} - 1+\Phi(\mu)\left(\bigcup\limits_{e\in E_c}\L_0[e]\right)=\sum\limits_{j\in c}\sum\limits_{e\in E,\ i(e)=j}\Phi(\mu)\left(\L_0[e]\right)-\frac{\epsilon}{2}\\
&=&\sum\limits_{j\in c}\sum\limits_{e\in E,\ i(e)=j}\int p_ed\mu-\frac{\epsilon}{2}=\sum\limits_{j\in c}\mu\left(K_j\right)-\frac{\epsilon}{2}.
\end{eqnarray*}
 This completes the proof of (iv).

(v) follows immediately from (iv).
\hfill$\Box$

\begin{Corollary}\label{usdmc} 
     Suppose $\mathcal{M}$  is consistent, contractive, uniformly continuous with $b<\infty$ and has a  dominating Markov chain such that $c<\infty$. Then  $P(\mathcal{M})$ is not empty.
\end{Corollary}
{\it Proof.}  
The assertion follow by Lemma \ref{prl} and Corollary \ref{eimc} (ii).
\hfill$\Box$

Finally, we remark that the study of a random dynamical system via an equivalent Markov system has another flexibility. As the following simple lemma shows, the problem of determining whether a Markov systems has an invariant measure can be reduced to that on a subsystem, which might be easier, as, obviously, any Borel probability measure on  $\bigcup_{i\in S}K_i$ can be uniquely identified with a member of $P(\bigcup_{i\in S}K_i)$. For example, as Examples \ref{dMse}, \ref{gnce}, \ref{mcae} and \ref{pex} show, a non-empty $\Omega$ sometimes contains a closed Markov subsystem. We cover this situation by the following simple corollary.

\begin{Definition}
We say that $\mathcal{M}$ contains a {\it Markov subsystem} iff there exists $S\subset N$ such that $(K_{i(e)}, w_e, p_e)_{e\in i^{-1}(S)}$ is a Markov system on $\bigcup_{i\in S}K_i$. We call the Markov subsystem {\it closed} iff $\bigcup_{i\in S}K_i$ is closed in $K$.
\end{Definition}

\begin{Lemma}\label{rp}
    Suppose $\mathcal{M}$ has a {\it Markov subsystem} $(K_{i(e)}, w_e, p_e)_{e\in i^{-1}(S)}$, for some $S\subset N$, which has an invariant $\mu\in P(\bigcup_{i\in S}K_i)$. Then $\mu\in P(\mathcal{M})$.
\end{Lemma}
{\it Proof.}  
  Let $B\subset K$ be Borel. Then
\begin{eqnarray*}
U^*\mu(B) &=& \int \sum\limits_{e\in E}p_e1_B\circ w_ed\mu =   \int\limits_{\bigcup\limits_{i\in S}K_i} \sum\limits_{e\in  i^{-1}(S)}p_e1_{B\cap\bigcup\limits_{i\in S}K_i}\circ w_ed\mu\\
&=&\mu\left(B\cap\bigcup_{i\in S}K_i\right)=\mu(B).
\end{eqnarray*}
\hfill$\Box$

\begin{Corollary}\label{sndct}
 Suppose $\mathcal{M}$ contains a closed Markov subsystem which satisfies the conditions of Corollary \ref{eimc} (ii). Then $P(\mathcal{M})$ is not empty.
\end{Corollary}
{\it Proof.}
   The assertion follows by Corollary \ref{eimc} (ii) and Lemma \ref{rp}.
 \hfill$\Box$

\section{Examples and applications}\label{easec}

In this section, in particular, some simple examples are given to which  the previous theory apparently  could not be applied.

\begin{Example}\label{gmc}
 This is to demonstrate that Proposition \ref{cml}, Corollary \ref{ndce} and Theorem \ref{est}  cover Theorem 2.1 in \cite{Wal}. 

Let $G:=(V,E,i,t)$ be a finite  directed graph. Set $\Sigma^-_G:=\{(...,\sigma_{-1},\sigma_0)|\  \sigma_m\in E\mbox{ and }  t(\sigma_m)=i(\sigma_{m-1})\ \mbox{ for all }
 m\in\mathbb{Z}\setminus\mathbb{N}\}$ (be the {\it  one-sided
  subshift of finite type} associated with $G$) endowed with the
  metric $d(\sigma,\sigma'):=2^k$ where $k$ is the smallest
  integer with $\sigma_i=\sigma'_i$ for all $k<i\leq 0$.  Let $T: \Sigma^-_G\longrightarrow \Sigma^-_G$ be the right shift map given by
 $(T\sigma)_i = \sigma_{i-1}$ for all $i\leq 0$. Let $g$
  be a positive continuous function on $\Sigma_G^-$ such that
  \[\sum\limits_{y\in T^{-1}(\{x\})}g(y)=1\mbox{ for all }x\in\Sigma_G^-.\]
  Set
  $K_i:=\left\{\sigma\in\Sigma^-_G:t(\sigma_0)=i\right\}$  for
  every $i\in V$ and, for
  every $e\in E$,
  \[ w_e(\sigma):=(...,\sigma_{-1},\sigma_{0},e),\  p_e(\sigma):=g(...,\sigma_{-1},\sigma_{0},e)
  \mbox{ for all }\sigma\in K_{i(e)}.\]
  Obviously, maps $( w_e|_{K_{i(e)}})_{e\in E}$ are contractions with a contraction rate $a=1/2$. Therefore,
  $\mathcal{M}_g:=\left(K_{i(e)}, w_e, p_e\right)_{e\in E}$
   defines a uniformly continuous contractive Markov system. Since each   $K_i$ is open, $\mathcal{M}_g$ is non-degenerate ($R1=0$). Therefore, by Corollary \ref{eimc} (i),  it has an invariant Borel probability measure. An invariant probability measure of $\mathcal{M}_g$ is called a $g$-{\it measure} \cite{Ke}. 
   Let $U_g$ be the Markov operator associated with  $\mathcal{M}_g$. Then, for every $f\in\mathcal{L}^B(\Sigma_G^-)$,
\[U_gf(x)=\sum\limits_{y\in T^{-1}(\{x\})}g(y)f(y)\mbox{ for all }x\in\Sigma_G^-.\]
Observe that,  in this case, $F$ is nothing else but the natural projection $\Sigma_G\longrightarrow\Sigma_G^-$ and $\Phi$ is the natural extension of a $g$-measure.  Moreover, in this example,  Corollary \ref{eimc} (v) is obvious, Theorem \ref{et} (iii) can be strengthened to globally  H\"{o}lder continuous $F$, and  Proposition \ref{cml}, Corollary \ref{ndce} (ii) and Theorem \ref{est} reduce to Theorem 2.1 in \cite{Wal} as follows. Let $\mathcal{B}$ denote the Borel $\sigma$-algebra on  $\Sigma_G^-$ and $ P_T(\Sigma_G^-)$ denote the set of all $T$-invariant members of $P(\Sigma_G^-)$.

\begin{Theorem}[Ledrappier, 1974 \cite{Le}]\label{Lt}
   Let $m\in P(\Sigma_G^-)$. Then the following are equivalent:\\
 (i) ${U_g}^*m=m$,\\
 (ii) $m\in P_T(\Sigma_G^-)$ and $E_m(f|T^{-1}\mathcal{B})= \sum_{z\in T^{-1}\{Tx\}} g(z)f(z)$ $m$-a.e. for all  $f\in\mathcal{L}^1(m)$,\\
 (iii) $m\in P_T(\Sigma_G^-)$ and $m$ is an equilibrium state for $\log g$.
 \end{Theorem}
{\it Proof.} Let  ${U_g}^*m=m$. Then, by Corollary \ref{ndce} (ii), $\Phi(m)\in E(\mathcal{M}_g)$, i.e. $\Phi(m)$ is $S$-invariant and 
 \begin{equation*}
      E_{\Phi(m)}\left(1_{_1[e]}|\mathcal{F}\right) =   p_e\circ F\ \ \ \Phi(m)\mbox{-a.e.} 
\end{equation*}
for all $e\in E$. Hence, $m$ is $T$-invariant, as $F\circ S^{-1} = T\circ F$. Let $e\in E$ and $A\in T^{-1}\mathcal{B}$. Let $_0[e]^-$ denote the cylinder set in $\Sigma_G^-$. Let $B\in\mathcal{B}$ such that $A = T^{-1}(B)$. Then
\begin{eqnarray}\label{eee}
  && \int\limits_{A}1_{_0[e]^-}dm= \int\limits_{F^{-1}(A)}1_{_0[e]}d\Phi(m)=\int\limits_{S^{-1}\left(F^{-1}(A)\right)}1_{_1[e]}d\Phi(m)\nonumber\\
&=&\int\limits_{F^{-1}(B)}p_e\circ F d\Phi(m)=\int\limits_{B}p_e dm=\int\limits_{B} \sum_{\sigma\in T^{-1}\{x\}} g(\sigma)1_{_0[e]^-}(\sigma) dm(x)\nonumber\\
&=&\int\limits_{A} \sum_{\sigma\in T^{-1}\{Tx\}} g(\sigma)1_{_0[e]^-}(\sigma) dm(x).
\end{eqnarray}
Hence, for any other cylinder set $_k[e_k,...,e_0]^-\subset\Sigma_G^-$,
\begin{eqnarray*}
   \int\limits_{A}1_{_k[e_k,...,e_0]^-}dm&=&  \int\limits_{A\cap _k[e_k,...,e_{-1}]^-}1_{_0[e_0]^-}dm\\
&=&\int\limits_{A\cap _k[e_k,...,e_{-1}]^-} \sum_{\sigma\in T^{-1}\{Tx\}} g(\sigma)1_{_0[e_0]^-}(\sigma) dm(x) \\
&=&\int\limits_{A} \sum_{\sigma\in T^{-1}\{Tx\}} g(\sigma)1_{ _k[e_k,...,e_{0}]^-}(\sigma) dm(x).
\end{eqnarray*}
By linearity, we obtain
\[  \int\limits_{A}sdm= \int\limits_{A} \sum_{\sigma\in T^{-1}\{Tx\}} g(\sigma)s(\sigma) dm(x)\]
for any simple function $s\in\mathcal{L}^1(m)$. Since the simple functions are dense in  $\mathcal{L}^1(m)$, we conclude that
 \begin{equation*}
     E_m(f|T^{-1}\mathcal{B})(x)= \sum_{z\in T^{-1}\{Tx\}} g(z)f(z)\ \ \ \mbox{ for } m\mbox{-a.e.  }x\in\Sigma_G^-\mbox{ and all } f\in\mathcal{L}^1(m).
\end{equation*}
This show the implication from (i) to (ii).

Form (ii) and \eqref{eee}, we obtain
 \begin{equation*}
      E_{\Phi(m)}\left(1_{_1[e]}|\mathcal{F}\right) = p_e\circ F\ \ \ \Phi(m)\mbox{-a.e.} 
\end{equation*}
for all $e\in E$.  Therefore, $\Phi(m)\in E(\mathcal{M}_g)$. Hence, by Theorem \ref{est}, $\Phi(m)$ is an equilibrium state for $u$. Observe that $h_m(T)=h_{\Phi(m)}(S)$, since $\Phi(m)$ is the natural extension of $m$. Therefore,
 \begin{equation}\label{esee}
      h_m(T)=\int u d\Phi(m) = \int \log p_{\sigma_1}\circ F(\sigma)d\Phi(m)(\sigma) = \int \log g dm.
\end{equation}
Thus $m$ is an equilibrium state for $ \log g$. This proves the implication from (ii) to (iii).

Finally, by \eqref{esee}, $\Phi(m)\in E(u)$ if $m$ is an equilibrium state for $ \log g$. Hence, by Theorem \ref{est} and  Proposition \ref{cml}, ${U_g}^*m = m$. This completes the proof.
\hfill$\Box$
\end{Example}

\begin{Example}\label{ieex}
Consider the following random dynamical system $D_R:=((\mathbb{R}, |.|), w_n, p_n)_{n\geq 3}$ where 
\begin{eqnarray*}
    w_n(x):=Z\sqrt{\log 2}\sqrt{\log n}x+1\mbox{ and }p_n(x):=\frac{1}{Z}\frac{1}{n(\log n)^2}
\end{eqnarray*}
for all $x\in\mathbb{R}$ and $n\geq 3$ where $Z$ is the suitable normalizing factor such that $\sum_{n\geq 3}p_n = 1$. Then a simple computation shows that
\[\sum\limits_{n\geq 3}p_n|w_n(x)-w_n(y)|\leq\sqrt{\log 2}\int\limits_2^\infty\frac{\sqrt{\log t}}{t(\log t)^2}dt|x-y|=\frac{1}{2}|x-y|\]
for all $x,y\in\mathbb{R}$, i.e. $D_R$ is contractive with a contraction rate $1/2$.  Also, for any choice of $x_0\in\mathbb{R}$, 
\[b=\sum\limits_{n\geq 3}p_n|w_n(x_0)-x_0|\leq\frac{1}{2}x_0+|1-x_0|.\] 
It is, obviously, non-degenerate. Thus, by Corollary \ref{eimc} (i) and Corollary \ref{ndce}, $D_R$ has a unique invariant Borel probability measure $\mu$. However, one easily checks that the Bernoulli measure $\Phi(\mu)$ has infinite entropy, i.e. $\{M\in E(\mathcal{M})|\ h_S(M)<\infty\}$ is empty. 
\end{Example}

\begin{Example}\label{prdm}
Consider the random dynamical system $([0,1], w_{e}, p_e)_{e\in\{0,1\}}$ where 
\begin{eqnarray*}
    w_0(x):=\frac{1}{2}x,&\ & w_1(x):=\frac{1}{2}+\frac{1}{2}x,\\
    p_0(x):=x,&\ &p_1(x):=1-x
\end{eqnarray*}
for all $x\in [0,1]$. The following Markov partition makes the random dynamical system to a uniformly continuous Markov system with strictly positive probability functions.
Set $K_1:=\{0\}$, $K_2:=(0,1)$, $K_3:=\{1\}$, 
\begin{eqnarray*}
  && w_a:=w_1|_{K_1}, w_b:=w_0|_{K_2}, w_c:=w_1|_{K_2}, w_d:=w_0|_{K_3}\\
  &&p_a:=p_1|_{K_1}, p_b:=p_0|_{K_2}, p_c:=p_1|_{K_2}, p_d:=p_0|_{K_3},
\end{eqnarray*}
and  $i:\{a,b,c,d\}\longrightarrow\{1,2,3\}$ by $i(a):=1$, $i(b):= 2$, $i(c):=2$, $i(d):=3$.  Obviously, $R1=1_{\{0\}\cup\{1\}}$ and 
\[R^21=\sum\limits_{e\in\{a,b,c,d\}}\partial p_e1_{\{0\}\cup\{1\}}\circ \bar w_e = 0.\]
Therefore, by Lemma  \ref {ndcl}, Markov system $(K_{i(e)},w_e, p_e)_{e\in\{a,b,c,d\}}$  satisfies the conditions of Corollary \ref{eimc} (i). One can choose also an infinite Markov partition, e.g. $K_1:=\{0\}$, $K_0:=\{1\}$,  $K_i:=(1-1/2^{i-2}, 1-1/2^{i-1}]$ for all $i\geq 2$. Then 
one easily sees (by drawing the directed graph of the Markov system) that
\[R1=1_{\{0\}}+1_{\{1-\frac{1}{2}\}}+1_{\{1-\frac{1}{4}\}}+...,\]
\[R^21=1_{\{0\}}+\frac{1}{2}1_{\{1-\frac{1}{2}\}}+\frac{1}{4}1_{\{1-\frac{1}{4}\}}+\frac{1}{8}1_{\{1-\frac{1}{8}\}}+...\]
and
\[R^31=\frac{1}{2}1_{\{0\}}+\frac{1}{8}1_{\{1-\frac{1}{2}\}}+\frac{1}{32}1_{\{1-\frac{1}{4}\}}+....\]
Since $R^31<1$, by Lemma  \ref {ndcl}, the resulting Markov system is non-degenerate, and therefore, also satisfies the conditions of Corollary \ref{eimc} (i) for any choice of $x_i\in K_i$ for all $i\in\mathbb{N}$. Furthermore, one easily checks that the dominating Markov chain has $c<\infty$, and therefore, by Lemma \ref{prl},  Condition \ref{dtc} is satisfied for all $x_0\in K$. (Note that the dominating Markov chain  has a positively recurrent communication class.)
\end{Example}

\begin{Example}\label{dMse}
Consider the random dynamical system $([0,1], w_{e}, p_e)_{e\in\{0,1\}}$ where $w_0$ and $w_1$ as in Example \ref{prdm}, but
\begin{eqnarray*}
    p_0(x):=1-x,\ p_1(x):=x
\end{eqnarray*}
for all $x\in [0,1]$. The random dynamical system has an equivalent, proper Markov system $\mathcal{M}:=(K_{i(e)},w_e, p_e)_{e\in\{a,b,c,d\}}$ where 
 $K_1:=\{0\}$, $K_2:=(0,1)$, $K_3:=\{1\}$, 
\begin{eqnarray*}
  && w_a:=w_0|_{K_1}, w_b:=w_0|_{K_2}, w_c:=w_1|_{K_2}, w_d:=w_1|_{K_3}\\
  &&p_a:=p_0|_{K_1}, p_b:=p_0|_{K_2}, p_c:=p_1|_{K_2}, p_d:=p_1|_{K_3},
\end{eqnarray*}
and $i:\{a,b,c,d\}\longrightarrow\{1,2,3\}$ by $i(a):=1$, $i(b):= 2$, $i(c):=2$, $i(d):=3$.  If one draws the directed graph associated with the Markov system, one can see that $(K_{i(e)},p_e,w_e)_{e\in i^{-1}(\{1,3\})}$ forms a Markov subsystem which has more than one invariant probability measures, and so the Markov system. The latter follows also by Corollary \ref{sndct}. Moreover, note that $\partial p_a = 0$, $\partial p_b = 1_{\{0\}}$, $\partial p_c = 1_{\{1\}}$ and $\partial p_d = 0$. Hence,
\[R1=1_{\{0\}\cup\{1\}},\]
 $R(R1)=R1$, and therefore, $\Omega=\{0\}\cup\{1\}$. Set $\sigma^0:=(...,b,b,b,...)$. Then, obviously, $\sigma^0\in T_2$, but $F(\sigma^0)=0$, and therefore, $\sigma^0\in F^{-1}(K_1)$. Hence $\sigma^0\notin G$. Set $\Lambda:=\delta_{\sigma^0}$. Then, by Lemma \ref{esfg},  $\Lambda\notin E(\mathcal{M})$, but, obviously,  $F(\Lambda)=\delta_0\in P(\mathcal{M})$. Let $A\in\mathcal{F}$, then
\[\int\limits_A1_{_1[b]}d\Lambda = \Lambda(A)=\int\limits_A1_{K_1}\circ F1_{T_2}d\Lambda=\int\limits_A\partial p_{b}\circ F1_{T_{i(b)}}d\Lambda.\]
Hence, $\Lambda\in E_\perp(\mathcal{M})$. Thus, by Theorem \ref{cndc}, $\mathcal{M}$ is degenerate. However, by Theorem \ref{ucct}, it is consistent. Also, it can be seen by Theorem \ref{imct} (ii), as obviously, $1_{\Omega}Rf = 1_{\Omega}Uf$ for all $f\in\mathcal{L}^B(K)$, or by Lemma \ref{uccl}. (Interestingly, $\Lambda(F^{-1}(F(G)))\geq \Lambda(F^{-1}(\{0\}))=1$.) Thus, the Markov  system satisfies the conditions of Corollary \ref{eimc} (ii). 

This example also can be used, in order to see  that Condition \ref{csc} is not necessary for the consistency of a Markov system. Set $K'_0:=\{0\}$, $K'_1:=\{1\}$, $K'_2:=\bigcup_{i\geq 2, even}(1-1/2^{i-2},1-1/2^{i-1}]$ and $K'_3:=\bigcup_{i\geq 3, odd}(1-1/2^{i-2},1-1/2^{i-1}]$. Then, the restrictions of the maps and the probability functions form a finite uniformly continuous Markov system. In this case,  $1\in\Omega$ also, but, obviously, $R1(1) = 2>U1(1)$, as $(\bar K'_2\setminus  K'_2)\cap(\bar  K'_3\setminus  K'_3)=\{1\}$.
\end{Example}

\begin{Example}\label{gnce}
Consider the random dynamical system $([0,1], w_{e}, p_e)_{e\in\{0,1\}}$ where $w_0$ and $w_1$  as in Example \ref{prdm}, but
\begin{eqnarray*}
    p_0(x):=1-x,\ p_1(x):=x
\end{eqnarray*}
for all $x\in (0,1]$ and $p_0(0) := p_1(0) := 1/2$. In this case, the random dynamical system has an equivalent, contractive, uniformly continuous Markov system $\mathcal{M}:=(K_{i(e)},w_e, p_e)_{e\in\{a,b,c,d,e\}}$ where 
 $K_1:=\{0\}$, $K_2:=(0,1)$, $K_3:=\{1\}$, 
\begin{eqnarray*}
  && w_a:=w_0|_{K_1}, w_b:=w_0|_{K_2}, w_c:=w_1|_{K_2}, w_d:=w_1|_{K_3}, w_e:=w_1|_{K_1}\\
  &&p_a:=p_0|_{K_1}, p_b:=p_0|_{K_2}, p_c:=p_1|_{K_2}, p_d:=p_1|_{K_3}, p_e:=p_1|_{K_1}
\end{eqnarray*}
and $i:\{a,b,c,d,e\}\longrightarrow\{1,2,3\}$ by $i(a):=1$, $i(b):= 2$, $i(c):=2$, $i(d):=3$, $i(e):=1$. As in Example \ref{dMse}, one easily checks that $R1=1_{\{0\}\cup\{1\}}$ and $R(R1)=R1$, and therefore, $\Omega=\{0\}\cup\{1\}$. Furthermore, the same way as in Example \ref{dMse}, one sees that the measure $\Lambda:=\delta_{(...,b,b,b,...)}\in E_\perp(\mathcal{M})$, but, obviously, $F(\Lambda)=\delta_0\notin P(\M)$. Thus, in this case, the Markov system is not consistent. However, by the symmetry of the system, $\delta_{(...,c,c,c,...)}\in E_\perp(\mathcal{M})$, and $F(\delta_{(...,c,c,c,...)})=\delta_1\in P(\M)$. This case is covered by Corollary \ref{sndct} with the closed Markov subsystem on $K_3$.
\end{Example}

\begin{Example}\label{nse}
Let $Blim$ be a Banach limit. Consider the following random dynamical system $(K, w_{e}, p_e)_{e\in\{0,1\}}$ where $K:=\{(x_1,x_2,...)|\ x_i\in [0,1]\mbox{ for all }i\in\mathbb{N}\}$ equipped with the supremum norm and
\begin{eqnarray*}
    w_0(x):=\frac{1}{2}x,&\ & w_1(x):=\left(\frac{1}{2}, \frac{1}{2},...\right)+\frac{1}{2}x,\\
    p_0(x):= Blim(x),&\ &p_1(x):=1-p_0(x)
\end{eqnarray*}
for all $x\in K$.
 Recall that the space $(K, \|.\|_\infty)$ is not separable.  Since $Blim$ is continuous, it is Borel measurable. Set $K_1:=\{x\in K|\ Blim(x) = 0\}$, $K_3:=\{x\in K|\ Blim(x) = 1\}$ and $K_2:=K\setminus(K_1\cap K_3)$. Then the random dynamical system is equivalent to the following Markov system. Set 
 \begin{eqnarray*}
  && w_a:=w_1|_{K_1}, w_b:=w_0|_{K_2}, w_c:=w_1|_{K_2}, w_d:=w_0|_{K_3}\\
  &&p_a:=p_1|_{K_1}, p_b:=p_0|_{K_2}, p_c:=p_1|_{K_2}, p_d:=p_0|_{K_3},
\end{eqnarray*}
and $i:\{a,b,c,d\}\longrightarrow\{1,2,3\}$ by $i(a):=1$, $i(b):= 2$, $i(c):=2$, $i(d):=3$. By the continuity of $Blim$, $K_1$ and $K_3$ are closed and $\bar K_2 = K$. Hence $R1 = 1_{K_1\cup K_2}$ and $R^21 = 0$. Therefore, by Lemma  \ref {ndcl}, Markov system $(K_{i(e)},w_e, p_e)_{e\in\{a,b,c,d\}}$ satisfies the assumptions of Corollary \ref{eimc}  (i).
\end{Example}

The following example is probably the most useful one.
\begin{Example}\label{mcae}
  Fix $n\in\mathbb{N}$ and $a_i\in [0, 1]$ for all $i\in\{0,...,2^n-1\}$. Consider the random dynamical system $([0,1], w_{e}, p_e)_{e\in \{0,1\}}$ where $w_0$ and $w_1$ as in  Example \ref{prdm},
     \[p_0(x):=\sum\limits_{i=0}^{2^n-1} a_i1_{Q_i}(x)+a_{2^n-1}1_{\{1\}}(x),\]
 for all $x\in [0,1]$, where $Q_i:=[i/2^n,  (i+1)/2^n)$ for all $0\leq i\leq 2^n-1$,  and $p_1:=1-p_0$. Now, set 
$K_i:=Q_i$ for all $0\leq i\leq 2^n-1$ and $K_{2^n}:=\{1\}$. Then $K_0,...,K_{2^n}$ obviously form a Markov partition for the random dynamical system. 

In order to construct the Markov system associated with it, set $E'_0:= \{i|\ a_i>0, 0\leq i\leq 2^n-1\}$ and $w'_i:=w_0|_{K_i}$ and $p'_i:=p_0|_{K_i}$ for all $i\in E'_0$. Set $w'_{2^n}:=w_0|_{\{1\}}$, $p'_{2^n}:=p_0|_{\{1\}}$ and $E_0:=E'_0\cup {2^n}$  if $a_{2^n-1}>0$. Otherwise,  $E_0:=E'_0$. Let  $E'_1:=\{i|\ a_i<1, 0\leq i\leq 2^n-1\}$ and set $w'_{-i}:=w_1|_{K_i}$ and $p'_{-i}:=p_1|_{K_i}$ for all $i\in E'_1$. Set $w'_{-2^n}:=w_1|_{\{1\}}$, $p'_{-2^n}:=p_1|_{\{1\}}$ and $E_1:=E'_1\cup \{2^n\}$  if $a_{2^n-1}<1$. Otherwise,  $E_1:=E'_1$. Finally, set $E:=E_0\cup E_1$ and $i:\ E\longrightarrow\{0,...,2^n\}$ by $i(e):=|e|$ for all $e\in E$. Then $(K_{i(e)},w'_e, p'_e)_{e\in E}$ is clearly a contractive, uniformly continuous Markov system which is equivalent to $([0,1], w_{e}, p_e)_{e\in \{0,1\}}$. Obviously,
\[R1=\sum\limits_{i=1}^{2^n-1}1_{\left\{\frac{i+1}{2^n}\right\}}.\]
Hence,
 \begin{eqnarray*}
 &&R^{n+1}1=R^n(R1)\\
&=&\sum\limits_{e_1,...,e_n}\partial p_{e_1}\partial p_{e_2}\circ \bar w_{e_1}\partial p_{e_n}\circ \bar w_{e_{n-1}}\circ ...\circ  \bar w_{e_1}\left( \sum\limits_{i=1}^{2^n-1}1_{\left\{\frac{i+1}{2^n}\right\}}\right)\circ \bar w_{e_{n}}\circ ...\circ  \bar w_{e_1}.
\end{eqnarray*}
Since for every $e_1,...,e_n\in E$ except for $i(e_1)=0$ and $i(e_1)=2^n-1$ and $i(e_1)=2^n$,
\[\bar w_{e_{n}}\circ ...\circ  \bar w_{e_1}\left(\left[\frac{i(e_1)}{2^n},\frac{i(e_1)+1}{2^n}\right]\right)\subset\left(\frac{i}{2^n},\frac{i+1}{2^n}\right)\]
for some $0\leq i\leq 2^n-1$, and $\partial p_e(0) = 0$  for all $e\in E$, and $\partial p_e(1) = 0$  for all $e\in E$ with $i(e) \neq 2^n-1$,
\[R^{n+1}1=b_{n+1}1_{\{1\}}\]
for some $0\leq b_{n+1}\leq 1$. If $b_{n+1}<1$, then  the Markov system  is non-degenerate, and therefore, has an invariant measure by Lemma  \ref {ndcl} and Corollary \ref{eimc} (i).
Otherwise,
\[R^{n+2}1=R\left(R^{n+1}1\right)=(1-a_{2^n-1})1_{\{1\}},\]
and therefore, the  Markov system  is non-degenerate if  $a_{2^n-1}>0$, and it is consistent otherwise. Hence, it has an invariant measure by Corollary \ref{usdmc} or  Corollary \ref{sndct}. 
\end{Example}

\begin{Example}\label{pex}
    Let $\mathcal{D}_r := (\mathbb{R}, w_e, p_e)_{e \in\{ 0,1\}}$ be the random dynamical system
    where  $w_0$ and $w_1$ as in  Example \ref{prdm},
    \begin{equation*}
        p_0(x) = \left\{\begin{array}{cc}
         a, & x\in\mathbb{Q} \\
         b, & x\in \mathbb{R}\setminus \mathbb{Q}
        \end{array}\right.
    \end{equation*}
    for some $0\leq a,b\leq 1$, where $\mathbb{Q}$ denotes the rational numbers, and $p_1 = 1 - p_0$. Set $K_0 := \mathbb{Q}$ and $K_1 := \mathbb{R}\setminus\mathbb{Q}$. Then clearly, $\{K_0,K_1\}$ is a Markov partition for $\mathcal{D}_r$, which makes it to a uniformly continuous Markov system. For this Markov system,  $\Omega=\mathbb{R}$. Therefore, Condition \ref{csc} is satisfied if and only if $Rf\leq Uf$ for all $f\in\mathcal{L}^B(\mathbb{R})$. However, one easily checks that $Rf=(a1_{K_1}+b1_{K_0})f\circ w_0+(1-(a1_{K_1}+b1_{K_0}))f\circ w_1$ for all $f\in\mathcal{L}^B(\mathbb{R})$. Hence, Condition \ref{csc} is satisfied if and only if $a=b$. Clearly, in this case, the Markov system has an invariant measure, in the agreement with Corollary \ref{usdmc}.
\end{Example}

\subsection*{Acknowledgements} 
The author would like to thank an anonymous reviewer for pointing out an error in the statement of a previous version of Corollary \ref{eimc} (iii), correcting numerous misprints and English grammar errors and suggestions on improvements to the text of this article.
Also, the author would like to thank Boris M. Gurevich for the invitations to give several talks on the subject at the Ergodic Theory and Statistical Mechanics Seminar at the Lomonosov Moscow State University and also the other organizers and participants of the seminar for the questions which helped to improve the article, in particular Sergey A. Pirogov whose question on the dominating Markov chain condition, the author strongly suspects, was an indirect way to point out that the formulation of the condition was unnecessarily strong.


\begin{thebibliography}{9}

\bibitem{BLLC} A. Baraviera, C. F. Lardizabal, A. O. Lopes, M. Terra Cunha,
A dynamical point of view of Quantum Information: entropy, pressure and Wigner measures,
in \emph{Dynamics, Games and Science II,
Springer Proceedings in Mathematics}  
\textbf{2} (2011) 161-185.

\bibitem{BDEG} M. F. Barnsley, S. G. Demko, J. H. Elton and J. S. Geronimo, Invariant measure
    for Markov processes arising from iterated function systems with place-dependent
    probabilities, 
\emph{ Ann. Inst. Henri Poincar\'{e}} 
\textbf{24}  (1988) 367-394.

\bibitem{BDEGE} M. F. Barnsley, S. G. Demko, J. H. Elton and J. S. Geronimo, Erratum: Invariant measure
    for Markov processes arising from iterated function systems with place-dependent
    probabilities, 
\emph{ Ann. Inst. Henri Poincar\'{e}}
 \textbf{25}  (1989) 589-590.

 \bibitem{Bog}   V. I. Bogachev, Measure theory. Vol. I,II.
 \emph{Springer} (2007).

\bibitem{Elton} J. H. Elton, An ergodic theorem for iterated maps, 
\emph{Ergod. Th. \& Dynam. Sys.} 
\textbf{7} (1987)  481-488.

\bibitem{DU} M. Denker and M.  Urba\'{n}ski, On the existence of conformal measures, 
\emph{Trans. Am. Math. Soc.}, 
\textbf{328} (1991) no. 2, 563-587.

\bibitem{FHK} H. F\"{o}llmer, U. Horst, A. Kirman, Equilibria in financial markets with heterogeneous agents: a probabilistic perspective, \emph{Journal of Mathematical Economics} \textbf{41} (2005) 123-155.

\bibitem{HS} K. Horbacz and T. Szarek, Irreducible Markov systems on Polish spaces,
\emph{Studia Math.}
\textbf{177} (2006), no. 3, 285--295.

\bibitem{Is} R. Isaac, Markov processes and unique stationary probability measures,
\emph{ Pacific J. Math.} 
\textbf{ 12} (1962)  273-286.

\bibitem{JOP}     A. Johansson, A. \"{O}berg, M. Pollicott, Countable state shifts and uniqueness of g-measures,
\emph{Amer. J. Math.}
\textbf{129} (6) (2007)  1501-1511.

\bibitem{Ke} M. Keane, Strongly Mixing $g$-Measures,
\emph{Inventiones math.} 
\textbf{16} (1972) 309-324.

\bibitem{Le} F. Ledrappier, Principe variationnel et syst\`{e}mes dynamiques symboliques, 
\emph{Z. Wahrscheinlichkeitstheorie verw. Gebiete} 
\textbf{30} (1974) 185-202.

\bibitem{MW} R.  Mauldin and S.  Williams, Hausdorff dimension in
graph directed constructions,\emph{Tran. AMS}
\textbf{309} (1988) 811-829.

\bibitem{MU0} R. Mauldin and  M. Urba\'{n}ski,  Gibbs states on the symbolic space over an infinite alphabet, 
\emph{Israel J. Math.}, 
\textbf{125} (December 2001) 93-130.

\bibitem{MU} R.  Mauldin and M. Urba\'{n}ski, 
\emph{Graph Directed Markov Systems: Geometry and Dynamics of Limit Sets}, Cambridge Tracts in Mathematics, volume 148, Cambridge University Press (2003).

\bibitem{MSU} V. Mayer, B. Skorulski, and M. Urba\'{n}ski, \emph{Distance expanding random mappings, thermodynamical formalism, Gibbs measures and fractal geometry}, Lecture Notes in Mathematics, volume 2036 (2011), Springer.

\bibitem{RU} M. Roy and M. Urba\'{n}ski, Random graph directed Markov systems,
\emph{Discrete Contin. Dyn. Syst.}, 
\textbf{30} (2011) no. 1 261-298.

\bibitem{Sa1}  O. Sarig, Thermodynamic formalism for countable Markov shifts,
\emph{Ergod. Th. \& Dynam. Sys.}  
\textbf{19} (1999) 1565-1593.

\bibitem{Sa2} O. Sarig,  Thermodynamic formalism for null recurrent potentials, \emph{Israel Journal of Mathematics}
 \textbf{121} no. 1 (2001) 285-311.

\bibitem{Sl} W. Slomczynski, Dynamical entropy, Markov operators, and iterated function
systems, 
\emph{Rozprawy Habilitacyjne Uniwersytetu Jagiello\'{n}skiego Nr
362, Wydawnictwo Uniwersytetu Jagiello\'{n}skiego}, Krak\'{o}w (2003).

 \bibitem{SU} B. Stratmann and M. Urba\'{n}ski, Pseudo-Markov systems and infinitely generated Schottky groups,
 \emph{Amer. J. Math.}, 
\textbf{129} (2007) n. 4, 1019-1062.

\bibitem{Sz} T. Szarek, Invariant measures for nonexpansive Markov operators on Polish
spaces, 
\emph{Diss. Math.} 
\textbf{415}, 1-62 (2003).

\bibitem{Wal} P. Walters, Ruelle's Operator Theorem and $g$-measures,
\emph{Tran. AMS} \textbf{214} (1975) 375-387.

\bibitem{Wer0} I. Werner, Invariant measures for some Markov processes arising from constructions of fractals (in German),
\emph{Diploma thesis at the University of Heidelberg} (January 6, 2002).

\bibitem{Wer1} I. Werner, Contractive Markov systems,
\emph{J. London Math. Soc. }
\textbf{71} (2005), no. 1, 236-258.

\bibitem{Wer3} I. Werner, Coding map for a contractive Markov system,
\emph{Math. Proc. Camb. Phil. Soc.}
\textbf{140} (2) (2006) 333-347,
\href{http://arxiv.org/abs/math/0504247}{ 	arXiv:math/0504247}.

\bibitem{Wer6} I. Werner, The generalized Markov measure as an equilibrium state,
\emph{Nonlinearity} 
\textbf{18} (2005) 2261-2274,
\href{http://arxiv.org/abs/math/0503644}{ 	arXiv:math/0503644}.

\bibitem{Wer10} I. Werner, Dynamically defined measures and equilibrium states,
\emph{J. Math. Phys.} 
\textbf{52} 122701 (2011),
\href{http://arxiv.org/abs/1101.2623}{ arXiv:1101.2623}.

\bibitem{Wer12} I. Werner, Erratum: Dynamically defined measures and equilibrium states,
\emph{J. Math. Phys.} 
\textbf{53}  079902 (2012),
\href{http://arxiv.org/abs/1101.2623}{ arXiv:1101.2623}.

\bibitem{Wer13} I. Werner, Contractive Markov systems II,
\href{http://arxiv.org/abs/math/0503633}{arXiv:math/0503633}.

\bibitem{Wer11} I. Werner, Fundamental Markov systems,
\href{http://arxiv.org/abs/math/0509120}{arXiv:math/0509120}.

\end{thebibliography}
\end{document}